\documentclass[review]{elsarticle}

\usepackage{lineno,hyperref}

\journal{Journal of \LaTeX\ Templates}

\begin{document}

\begin{frontmatter}
\title{A critical branching process with immigration in random environment\tnoteref{mytitlenote}}
\tnotetext[mytitlenote]{This work is supported by the Russian Science Foundation under grant 19-11-00111.}
\author{Afanasyev V.I.}
\address{Steklov Mathematical Institute of Russian Academy of Sciences, 8 Gubkin St., Moscow, 119991, Russia Email: viafan@mi.ras.ru}
\begin{abstract}
A Galton-Watson branching process with immigration evolving in a random environment is considered. Its associated random walk is assumed to be oscillating. We prove a functional limit theorem in which the process under consideration is normalized by a random coefficient depending on the random environment only. The distribution of the limiting process is described in terms of a strictly stable Levy process and a sequence of independent and identically distributed random variables which is independent of this process.
\end{abstract}
\end{frontmatter}


\textbf{1}. \textbf{Introduction and statement of main result\medskip }

Let $\left( \Omega ,\mathcal{F},\mathbf{P}\right) $ be a probability space
and $\Delta $ be the space of probability measures on $\mathbf{N}%
_{0}:=\left\{ 0,1,\ldots \right\} $ equipped with the metric of total
variation. \textit{A} \textit{random environment} is a sequence of random
elements $Q_{1},Q_{2},\ldots $, mapping the space $\left( \Omega ,\mathcal{F}%
,\mathbf{P}\right) $ into $\Delta ^{2}$. Thus, $Q_{n}$ for each $n\in
\mathbf{N}$ has the form $\left( F_{n},G_{n}\right) $, where $F_{n},G_{n}$
are probability measures on $\mathbf{N}_{0}$. \textit{A branching process
with immigration in random environment} ((BPIRE)) is a stochastic process
possessing the following properties. For a fixed random environment $\left\{
Q_{n}, n\in \mathbf{N}\right\} $ this is an inhomogeneous branching
Galton-Watson process with immigration (see \cite{Athreya1972}, Chapter 6,
\S\ 7). Here, for each $n\in \mathbf{N}$, the number of immigrants joining
the $\left( n-1\right) $th generation has the distribution $G_{n}$ and the
offspring reproduction law of particles of the $\left( n-1\right) $th
generation is $F_{n}$.

Let $Z_{n}$ be the size of $n$th generation without the immigrants which
joined this generation (we assume that $Z_{0}=0$), $\eta _{n}$ be the number
of immigrants which joined the $n$th generation. Let $f_{n}\left( \cdot
\right) \ $and $g_{n}\left( \cdot \right) $ be generating functions of
distributions $F_{n}\ $and $G_{n}$ respectively.

We consider this model under the assumption that the random elements $%
Q_{1},Q_{2},\ldots $ are independent and identically distributed. A more
detailed definition of the BPIRE can be found in \cite{Dyakonova}.

Set for $i\in \mathbf{N}$
\[
X_{i}=\ln f_{i}^{\prime }\left( 1\right) ,\qquad \mu
_{i}=g_{i}^{\prime }\left( 1\right)
\]%
(suppose that $0<f_{1}^{\prime }\left( 1\right) <+\infty $, $0<g_{i}^{\prime
}\left( 1\right) <+\infty $ a.s.). Introduce the so-called \textit{%
associated random walk}:
\[
S_{0}=0,\qquad S_{n}=\sum\limits_{i=1}^{n}X_{i},\, n\in \mathbf{N%
}.
\]%
It is clear that the random vectors $\left( X_{1},\mu _{1}\right) ,\left(
X_{2},\mu _{2}\right) ,\ldots $ are independent and identically distributed
under our assumptions.

We impose the following restriction on the distribution of $X_{1}$.

\textbf{Hypothesis A}\textit{.} The distribution of $X_{1}$ belongs without
centering to the domain of attraction of some stable law with index $\alpha
\in \left( 0,2\right] $ and the limit law is not a one-sided stable law.

Under Hypothesis A the Skorokhod functional limit theorem is valid (see, for
instance, \cite{Kallenberg2002}, Chapter 16): there are such positive
normalizing constants $C_{n}$ that, as $n\rightarrow \infty $,%
\begin{equation}
W_{n} \stackrel{D}{\to}W,  \label{1}
\end{equation}%
where $W_{n}=\left\{ C_{n}^{-1}S_{\left\lfloor nt\right\rfloor },\,
t\geq 0\right\} $, the process $W=\left\{ W\left( t\right) ,\, t\geq
0\right\} $ is a strictly stable Levy process with index $\alpha \in
\left( 0,2\right] $ and the symbol $\stackrel{D}{\to}$ means
convergence in distribution in the space $D\left[ 0,+\infty \right)
$ with Skorokhod topology. Moreover,
\[
C_{n}=n^{1/\alpha }l\left( n\right) ,
\]%
where $\left\{ l\left( n\right), n\in \mathbf{N}\right\} $ is a
slowly varying sequence. It is known that the finite-dimensional
distributions of the process $W$ are absolutely continuous. Note that $\rho
:=\mathbf{P}\left( W\left( 1\right) >0\right) \in \left( 0,1\right) $ given
Hypothesis A. Thus, the Spitzer-Doney condition is satisfied:%
\begin{equation}
\lim_{n\rightarrow \infty }\mathbf{P}\left( S_{n}>0\right) =\rho \in
\left( 0,1\right) . \label{2}
\end{equation}

The Spitzer-Doney condition means that the random walk $\left\{
S_{n}\right\} $ is oscillating. As result, the absolute values of its strict
descending ladder heights constitute a renewal process with the
corresponding renewal function $v\left( x\right) $, $x\geq 0$ (see \cite%
{Afanasyev2005645} for a detailed definition of the function $v\left( \cdot
\right) $). Similarly, weak ascending ladder heights of the random walk $%
\left\{ S_{n}\right\} $ generate a renewal process with the corresponding
renewal function $u\left( x\right) $, $x\geq 0$.

The aim of this paper is to prove a functional limit theorem for the process
$\left\{ Z_{\left\lfloor nt\right\rfloor },\, t\geq 0\right\} $, as $%
n\rightarrow \infty $ (see Theorem 1).

We need some notation and definitions to formulate the theorem. Let for $%
n\in \mathbf{N}$%
\[
M_{n}=\max_{1\leq i\leq n}S_{i},\qquad L_{n}=\min_{%
0\leq i\leq n}S_{i}.
\]%
It is known (see, for instance, \cite{Afanasyev2005645}, Lemma 2.5) that, if
the Spitzer-Doney condition (2) is satisfied, then, as $n\rightarrow \infty $%
,%
\begin{equation}
\left\{ \left. \left( Q_{i},S_{i},\mu _{i}\right) ,\, i\in \mathbf{N}%
\, \right\vert \, L_{n}\geq 0\right\} \stackrel{D}{\to}%
\left\{ \left( Q_{i}^{+},S_{i}^{+},\mu _{i}^{+}\right) ,\,%
i\in \mathbf{N}\right\} ,  \label{3}
\end{equation}%
\begin{equation}
\left\{ \left. \left( Q_{i},S_{i},\mu _{i}\right) ,\,i\in \mathbf{N}%
\right\vert \,M_{n}\,<0\right\} \stackrel{D}{\to}\left\{
\left( Q_{i}^{-},S_{i}^{-},\mu _{i}^{-}\right) ,\,i\in \mathbf{N}%
\right\} ,  \label{4}
\end{equation}%
where $\left\{ \left( Q_{i}^{+},S_{i}^{+},\mu _{i}^{+}\right) \right\} $, $%
\left\{ \left( Q_{i}^{-},S_{i}^{-},\mu _{i}^{-}\right) \right\} $ are some
random sequences. Moreover: a) the sequences $\left\{ Q_{i}^{+},\,i\in
\mathbf{N}\right\} $, $\left\{ Q_{i}^{-},\,i\in \mathbf{N}\right\} $
can be viewed as some random environments; b) the sequences $\left\{
S_{i}^{+},\,i\in \mathbf{N}\right\} $, $\left\{ S_{i}^{-},\,i\in
\mathbf{N}\right\} $ are the corresponding associated random walks ($%
S_{0}^{+}=S_{0}^{-}=0$); c) the sequences $\left\{ \mu _{i}^{+},\,i\in
\mathbf{N}\right\} $ and $\left\{ \mu _{i}^{-},\,i\in \mathbf{N}%
\right\} $ are positive and constructed by $\left\{ Q_{i}^{+},\,i\in
\mathbf{N}\right\} $ and $\left\{ Q_{i}^{-},\,i\in \mathbf{N}\right\} $%
, respectively, the same as the sequence $\left\{ \mu _{i},\,i\in
\mathbf{N}\right\} $ is constructed by $\left\{ Q_{i},\,i\in \mathbf{N}%
\right\} $. Suppose that the sequences $\left\{ Q_{i}^{+},\,i\in
\mathbf{N}\right\} $, $\left\{ Q_{i}^{-},\,i\in \mathbf{N}\right\} $
are defined on the same probability space $\left( \Omega ^{\ast },\mathcal{F}%
^{\ast },\mathbf{P}^{\ast }\right) $ and are independent (below we denote
the expectation on this probability space by $\mathbf{E}^{\ast }$).

We now come back to our initial BPIRE. Set $\mathbf{N}_{i}=\left\{
i,i+1,\ldots \right\} $ for $i\in \mathbf{Z}$. Fix $i\in \mathbf{N}_{0}$
and, for $n\in \mathbf{N}_{i}$, denote by $Z_{i,n}$ the total number of
particles in the $n$th generation which are the descendants of the
immigrants joined the $i$th generation (we assume that $Z_{i,n}=0$ for $%
i\geq n$ and $i<0$). Note that the random sequence $\left\{ \eta _{i};\, Z_{i,n},\,n\in \mathbf{N}_{i+1}\right\} $ is a usual (without
immigration) branching process in the random environment $\left\{ G_{i+1};%
\,F_{n},\,n\in \mathbf{N}_{i+1}\right\} $. In particular, if the
random environment is fixed, then $G_{i+1}$ is the distribution of the
random variable $\eta _{i}$ which should be interpreted as the number of
particles in the initial generation. Set for $n\in \mathbf{N}_{i}$
\[
a_{i,n}=e^{-\left( S_{n}-S_{i}\right) }.
\]
The sequence $\left\{ \eta _{i};\,a_{i,n}Z_{i,n},\,n\in \mathbf{N%
}_{i+1}\right\} $ is a nonnegative martingale if the random environment $%
\left\{ G_{i+1};\,F_{n},\,n\in \mathbf{N}_{i+1}\right\} $ is
fixed. Hence (without assuming that the random environment is fixed), there
is a finite limit $\lim_{n\rightarrow \infty }a_{i,n}Z_{i,n}$ $\mathbf{P}$%
-a.s.

Set%
\[
Q_{i}^{\ast }=\left\{
\begin{array}{c}
Q_{i}^{+},\qquad i\in \mathbf{N}, \\
Q_{-i+1}^{-},\qquad i\in \mathbf{Z\setminus N,}%
\end{array}%
\right.
\]%
\[
S_{i}^{\ast }=\left\{
\begin{array}{c}
S_{i}^{+},\qquad\qquad i\in \mathbf{N}_{0}\mathbf{,} \\
-S_{-i}^{-},\qquad i\in \mathbf{Z\setminus N}_{0}\mathbf{,}%
\end{array}%
\right.
\]%
\[
\mu _{i}^{\ast }=\left\{
\begin{array}{c}
\mu _{i}^{+},\qquad\qquad i\in \mathbf{N}, \\
\mu _{-i+1}^{-},\qquad i\in \mathbf{Z\setminus N.}%
\end{array}%
\right.
\]%
The sequence $\mathcal{E}^{\ast }:=\left\{ Q_{k}^{\ast },\,k\in
\mathbf{Z}\right\} $ can be considered as a random environment (we denote
the components of $Q_{k}^{\ast }$ by $G_{k}^{\ast }\ $and $F_{k}^{\ast }$).
We assume that the probability space $\left( \Omega ^{\ast },\mathcal{F}%
^{\ast },\mathbf{P}^{\ast }\right) $ is reach enough for we are able to
define on it a branching process with immigration in the random environment $%
\mathcal{E}^{\ast }$. Fix $i\in \mathbf{Z}$ and, for $j\in \mathbf{N}_{i}$,
denote by $Z_{i,j}^{\ast }$ the total number of particles in the $j$th
generation being descendants of immigrants which joined the $i$th generation
(we denote the number of such immigrants as $\eta _{i}^{\ast }$). Note that
the sequence $\left\{ \eta _{i}^{\ast };\,Z_{i,j}^{\ast },\,j\in
\mathbf{N}_{i+1}\right\} $ is a branching process in the random environment $%
\left\{ G_{i+1}^{\ast };\,F_{j}^{\ast },\,j\in \mathbf{N}%
_{i+1}\right\} $ with the initial value $\eta _{i}^{\ast }$. The sequence $%
\left\{ S_{j}^{\ast }-S_{i}^{\ast },\,j\in \mathbf{N}_{i}\right\} $ is
the associated random walk and the random variable $\mu _{i}^{\ast }$ is
under fixed environment the mean of the random variable $\eta _{i}^{\ast }$.
Set%
\[
a_{i,j}^{\ast }=e^{-\left( S_{j}^{\ast }-S_{i}^{\ast }\right) }.
\]%
In accordance with the above the limit%
\begin{equation}
\lim_{j\rightarrow \infty }a_{i,j}^{\ast }Z_{i,j}^{\ast }=:\zeta _{i}^{\ast }
\label{5}
\end{equation}%
exists $\mathbf{P}^{\ast }$-a.s. and $\mathbf{P}^{\ast }\left( \zeta
_{i}^{\ast }>0\right) >0$ for $i\in \mathbf{N}_{0}$ (see \cite%
{Afanasyev2005645}, Proposition 3.1).

Introduce the following random series:%
\[
\Sigma _{1}:=\sum\limits_{i\in \mathbf{Z}}\mu _{i+1}^{\ast }e^{-S_{i}^{\ast
}},\qquad \Sigma _{2}:=\sum\limits_{i\in \mathbf{Z}}\zeta _{i}^{\ast
}e^{-S_{i}^{\ast }}
\]%
It is clear that $\Sigma _{1}>0$ $\mathbf{P}^{\ast }$-a.s. and $\mathbf{P}%
^{\ast }\left( \Sigma _{2}>0\right) >0$. Both series converge $\mathbf{P}%
^{\ast }$-a.s. under certain restrictions (see Lemma 4).

Let $W$ be a strictly stable Levy process with index $\alpha $ (in
the sequel we call $W$ simply \textit{the Levy process}). By the
Levy process we
specify \textit{the (lower) level }$L=\left\{ L\left( t\right) ,\,%
t\geq 0\right\} $ \textit{of the Levy process} as%
\[
L\left( t\right) =\inf_{s\in \left[ 0,t\right] }W\left( s\right) .
\]%
Let, further, $\gamma _{1},\gamma _{2},\ldots $ be an independent of $W$
sequence of independent random variables distributed as the random variable $%
\Sigma _{2}/\Sigma _{1}$.

By these ingredients we define finite-dimensional distributions of a random
process $Y=\left\{ Y(t),\,t\geq 0\right\} $ which plays an important
role in the sequel. First we set $Y(0)=0$. Consider an arbitrary $m\in
\mathbf{N}$ and arbitrary moments $t_{1},t_{2},\ldots ,t_{m}$: $%
0=t_{0}<t_{1}<t_{2}<\ldots <t_{m}$. The random vector $\left\{
Y(t_{1}),\ldots ,Y(t_{m})\right\} $ coincides in distribution with the
following vector $\widehat{Y}:=\left\{ \widehat{Y}_{1},\ldots ,\widehat{Y}%
_{m}\right\} $. We describe at first the possible values of the vector $%
\widehat{Y}$. Its first several coordinates coincide with $\gamma _{1}$, the
next several coordinates coincide with $\gamma _{2}$ and so on up to the $m$%
th coordinate. The coordinates of the vector $\widehat{Y}$ are specified
according to the level $L$ of the Levy process $W$. The first coordinate $%
\widehat{Y}_{1}$ is equal to $\gamma _{1}$. Let the coordinate $\widehat{Y}%
_{k}$ for some $k<m$ be known. For instance, $\widehat{Y}_{k}=\gamma _{l}$
for some $l\in \mathbf{N}$. If the level of the Levy process at the moment $%
t_{k+1}$ remains the same as at moment $t_{k}$, i.e. $L\left( t_{k+1}\right)
=L\left( t_{k}\right) $, then $\widehat{Y}_{k+1}=\gamma _{l}$. If the level
of the Levy process at the moment $t_{k+1}$ is changed, i.e. $L\left(
t_{k+1}\right) <L\left( t_{k}\right) $, then $\widehat{Y}_{k+1}=\gamma
_{l+1} $.

Set for $n\in \mathbf{N}_{0}$
\[
a_{n}=e^{-S_{n}},\qquad b_{n}=\sum\limits_{i=0}^{n-1}\mu
_{i+1}e^{-S_{i}}\ (b_{0}=0).
\]%
Introduce for each $n\in \mathbf{N}$ the random process $Y_{n}=\left\{
Y_{n}\left( t\right) ,\,t\geq 0\right\} $, where
\[
Y_{n}\left( 0\right) =0,\qquad Y_{n}\left( t\right) =\frac{%
a_{\left\lfloor nt\right\rfloor }}{b_{\left\lfloor nt\right\rfloor }}%
Z_{\left\lfloor nt\right\rfloor }.
\]%
Note that for $k\in \mathbf{N}$ the ratio $b_{k}/a_{k}$ is equal to the mean
of $Z_{k}$ for a fixed random environment.

Let the symbol\textit{\ }$\Rightarrow $\textit{\ }means\textit{\ }%
convergence of random processes in the sense of finite-dimensional
distributions and $\ln ^{+}x=\max \left( 0,\ln x\right) $ for $x>0$.\medskip

\textbf{Theorem 1. }\textit{If Hypothesis A is valid and }$\mathbf{E}\left(
\ln ^{+}\mu _{1}\right) ^{\alpha +\varepsilon }<+\infty $ \textit{for some }$%
\varepsilon >0$,\textit{\ then},\textit{\ as }$n\rightarrow \infty $,%
\[
Y_{n}\Rightarrow Y.
\]%
\medskip

A detailed description of the theory of critical (when Hypothesis A is
valid) branching processes in random environment is available in \cite%
{Afanasyev2005645} and \cite{Kersting2017}.

A particular case of a subcritical BPIRE (when the offspring generating
function $f_{n}\left( \cdot \right) $ is fractional-linear and $g_{n}\left(
s\right) \equiv s$ for each $n\in \mathbf{N}$) was considered in \cite%
{Kesten1975145}. The main attention there was paid to obtaining an
exponential estimate for the tail distribution of the so-called life period
of this process (i.e., the time until the first extinction). A more general
class of subcritical BPIRE was analyzed in \cite{Key1987344} where a limit
theorem describing the population size at a distant moment was proved and an
exponential estimate for the tail distribution of the life period was
established. A strong law of large numbers and a central limit theorem for a
wide class of subcritical BPIRE were proved in \cite{Roitershtein20071573}.

A critical BPIRE was considered in \cite{Bauernschubert2014683} where
sufficient conditions of transience and recurrence were obtained. The author
of \cite{Afanasyev2013547}, studying a random walk in random environment,
proved a particular case of Theorem 1 (when the offspring generating
function $f_{n}\left( \cdot \right) $ is fractional-linear and $g_{n}\left(
s\right) \equiv s$ for each $n\in \mathbf{N}$). We would like to stress that
the proof used in the present paper differs significantly from that one
given in \cite{Afanasyev2013547}. We also mention the papers \cite%
{Afanasyev2014525}, \cite{Afanasyev2017178} and \cite{Afanasyev2018139} in
which critical and supercritical processes (with stopped immigration) are
considered under some restrictions on their lifetime.

Recent papers \cite{Dyakonova} and \cite{Li2020} contain exact asymptotic
formulae for the tail distribution of the life period for critical and
subcritical BPIRE.\medskip

\textbf{2.} \textbf{Auxiliary statements}\medskip

Let $\tau _{n}$ be the first moment when the minimum of the random walk $%
S_{0},\ldots ,S_{n}$ is attained:

\[
\tau _{n}=\min \left\{ i:S_{i}=L_{n},\,0\leq i\leq n\right\} .
\]%
Set for $n\in \mathbf{N}$%
\[
S_{i,n}^{\prime }=\left\{
\begin{array}{c}
S_{\tau _{n}+i}-S_{\tau _{n}},\qquad i\in \mathbf{N}_{\left( -\tau
_{n}\right) }, \\
0,\qquad\qquad i\in \mathbf{Z\setminus N}_{\left( -\tau
_{n}\right) }.%
\end{array}%
\right.
\]%
For positive integers numbers $n_{1}<n_{2}$ set
\[
L_{n_{1},n_{2}}=\min_{n_{1}\leq i\leq n_{2}}S_{i}.
\]%
\medskip

\textbf{Lemma 1}. \textit{If the Spitzer-Doney condition} (2) \textit{is
satisfied}, \textit{then}, \textit{as} $n\rightarrow \infty $,
\begin{equation}
\left\{ S_{i,n}^{\prime },\,i\in
\mathbf{Z}\right\}\stackrel{D}{\to}\left\{ S_{i}^{\ast },\,i\in
\mathbf{Z}\right\} . \label{6}
\end{equation}%
\textit{\medskip }

\textit{Proof}. We demonstrate for simplicity only convergence of
one-dimensional distributions. Fix $i\in \mathbf{N}_{0}$. Let $A$ be a
one-dimensional $S_{i}^{\ast }$-continuous (relative to the measure $\mathbf{%
P}^{\ast }$) Borel set. Then for $n\geq i$%
\begin{eqnarray*}
&&\mathbf{P}\left( S_{i,n}^{\prime }\in A,\,\tau _{n}+i\leq n\right) \\
&=&\sum\limits_{k=0}^{n-i}\mathbf{P}\left( S_{i,n}^{\prime }\in A,\,%
\tau _{n}=k\right) \\
&=&\sum\limits_{k=0}^{n-i}\mathbf{P}\left( S_{k+i}-S_{k}\in A,\,%
S_{k}<L_{k-1},\,S_{k}\leq L_{k+1,n}\right)
\end{eqnarray*}%
and by the Markov property of random walks we have that%
\begin{eqnarray*}
&&\mathbf{P}\left( S_{k+i}-S_{k}\in A,\,S_{k}<L_{k-1},\,%
S_{k}\leq L_{k+1,n}\right) \\
&=&\mathbf{P}\left( S_{k}<L_{k-1}\right) \mathbf{P}\left( S_{k+i}-S_{k}\in A,%
\,S_{k}\leq L_{k+1,n}\right) \\
&=&\mathbf{P}\left( S_{k}<L_{k-1}\right) \mathbf{P}\left( S_{i}\in A,\,%
L_{n-k}\,\geq 0\right) \\
&=&\mathbf{P}\left( \left. S_{i}\in A\,\right\vert \,L_{n-k}%
\,\geq 0\right) \mathbf{P}\left( S_{k}<L_{k-1}\right) \mathbf{P}\left(
L_{n-k}\,\geq 0\right) \\
&=&\mathbf{P}\left( \left. S_{i}\in A\,\right\vert \,L_{n-k}%
\,\geq 0\right) \mathbf{P}\left( \tau _{n}=k\right) .
\end{eqnarray*}

Thus,
\begin{equation}
\mathbf{P}\left( S_{i,n}^{\prime }\in A,\,\tau _{n}+i\leq n\right)
=\sum\limits_{k=0}^{n-i}\mathbf{P}\left( \left. S_{i}\in A\,%
\right\vert \,L_{n-k}\,\geq 0\right) \mathbf{P}\left( \tau
_{n}=k\right) .  \label{7}
\end{equation}

If the Spitzer-Doney condition is satisfied, then the following generalized
arcsine law is valid (see, for instance, \cite{Bingham1987}, Chapter 8,
Theorem 8.9.9): for $x\in \left[ 0,1\right] $
\begin{equation}
\lim_{n\rightarrow \infty }\mathbf{P}\left( \frac{\tau _{n}}{n}%
\leq x\right) =\frac{\sin \left( \pi \rho \right) }{\pi }\int%
\limits_{0}^{x}u^{\rho -1}\left( 1-u\right) ^{-\rho }du.  \label{8}
\end{equation}

We pass to the limit in formula (7), as $n\rightarrow \infty $. Due to (8)

\[
\lim_{n\rightarrow \infty}\mathbf{P}\left( \tau _{n}+i\leq n\right)
={\lim_{n\rightarrow \infty } }\mathbf{P}\left( \tau _{n}/n\leq
1-i/n\right) =1.
\]%
Therefore the limit of the left-hand side of (7) coincides with the limit of
probability $\mathbf{P}\left( S_{i,n}^{\prime }\in A\right) $, as $%
n\rightarrow \infty $, if at least one of these limits exists.

If $\varepsilon \in \left( 0,1\right) $ and $n$ is large enough, then by (7)
\begin{equation}
\mathbf{P}\left( S_{i,n}^{\prime }\in A,\,\tau _{n}+i\leq n\right)
=P_{1}\left( n,\varepsilon \right) +P_{2}\left( n,\varepsilon \right) ,
\label{9}
\end{equation}%
where%
\[
P_{1}\left( n,\varepsilon \right) =\sum\limits_{k=0}^{\left\lfloor \left(
1-\varepsilon \right) n\right\rfloor }\mathbf{P}\left( \left. S_{i}\in A%
\,\right\vert \,L_{n-k}\,\geq 0\right) \mathbf{P}\left(
\tau _{n}=k\right) ,
\]%
\[
P_{2}\left( n,\varepsilon \right) =\sum\limits_{k=\left\lfloor \left(
1-\varepsilon \right) n\right\rfloor +1}^{n-i}\mathbf{P}\left( \left.
S_{i}\in A\,\right\vert \,L_{n-k}\,\geq 0\right) \mathbf{P}%
\left( \tau _{n}=k\right) .
\]

Clearly,%
\begin{eqnarray*}
P_{2}\left( n,\varepsilon \right) &\leq &\sum\limits_{k=\left\lfloor \left(
1-\varepsilon \right) n\right\rfloor +1}^{n}\mathbf{P}\left( \tau
_{n}=k\right) =\mathbf{P}\left( \tau _{n}>\left\lfloor \left( 1-\varepsilon
\right) n\right\rfloor \right) \\
&&\stackrel{n\rightarrow \infty}{\longrightarrow}1-\frac{\sin \left(
\pi \rho \right) }{\pi }\int\limits_{0}^{1-\varepsilon }u^{\rho
-1}\left( 1-u\right) ^{-\rho }du\stackrel{\varepsilon\rightarrow
0}{\longrightarrow}0.
\end{eqnarray*}%
Therefore%
\begin{equation}
\lim_{\varepsilon \rightarrow 0}{%
\limsup_{n\rightarrow \infty } }P_{2}\left( n,\varepsilon \right)
=0. \label{10}
\end{equation}

In view of (3) the probability $\mathbf{P}\left( \left. S_{i}\in A\,%
\right\vert \,L_{n-k}\,\geq 0\right) $ tends, as $n\rightarrow
\infty $, to $\mathbf{P}^{\ast }\left( S_{i}^{+}\in A\,\right) =%
\mathbf{P}^{\ast }\left( S_{i}^{\ast }\in A\,\right) $ uniformly over $%
0\leq k\leq \left\lfloor \left( 1-\varepsilon \right) n\right\rfloor $.
Consequently,%
\begin{eqnarray*}
\lim_{n\rightarrow \infty }P_{1}\left( n,\varepsilon \right) &=&%
\mathbf{P}^{\ast }\left( S_{i}^{\ast }\in A\,\right) {\lim_{%
n\rightarrow \infty } }\sum\limits_{k=0}^{\left\lfloor \left(
1-\varepsilon \right) n\right\rfloor }\mathbf{P}\left( \tau _{n}=k\right) \\
&=&\mathbf{P}^{\ast }\left( S_{i}^{\ast }\in A\,\right) {\lim_{%
n\rightarrow \infty } }\mathbf{P}\left( \tau _{n}\leq \left\lfloor
\left( 1-\varepsilon \right) n\right\rfloor \right) \\
&=&\mathbf{P}^{\ast }\left( S_{i}^{\ast }\in A\,\right) \frac{\sin
\left( \pi \rho \right) }{\pi }\int\limits_{0}^{1-\varepsilon }u^{\rho
-1}\left( 1-u\right) ^{-\rho }du \\
&&\stackrel{\varepsilon\rightarrow
0}{\longrightarrow}\mathbf{P}^{\ast }\left( S_{i}^{\ast }\in
A\,\right)
\end{eqnarray*}%
implying%
\begin{equation}
{\lim_{\varepsilon \rightarrow 0} }{%
\lim_{n\rightarrow \infty } }P_{1}\left( n,\varepsilon \right)
=\mathbf{P}^{\ast }\left( S_{i}^{\ast }\in A\,\right) . \label{11}
\end{equation}

It follows from relations (9)-(11) that for $i\in \mathbf{N}_{0}$

\[
{\lim_{n\rightarrow \infty } }\mathbf{P}\left( S_{i,n}^{\prime }\in
A,\,\tau _{n}+i\leq n\right) =\mathbf{P}^{\ast }\left( S_{i}^{\ast
}\in A\,\right) .
\]%
Thus,

\begin{equation}
{\lim_{n\rightarrow \infty }}\mathbf{P}\left( S_{i,n}^{\prime }\in
A\right) =\mathbf{P}^{\ast }\left( S_{i}^{\ast }\in A\,\right) .
\label{12}
\end{equation}

We now fix $i\in \mathbf{N}$. Let $A$ be a one-dimensional $S_{-i}^{\ast }$%
-continuous (relative to the measure $\mathbf{P}^{\ast }$) Borel set. Then
for $n\geq i$%
\begin{eqnarray*}
&&\mathbf{P}\left( S_{-i,n}^{\prime }\in A,\,\tau _{n}-i\geq 0\right)
\\
&=&\sum\limits_{k=i}^{n}\mathbf{P}\left( S_{-i,n}^{\prime }\in A,\,%
\tau _{n}=k\right) \\
&=&\sum\limits_{k=i}^{n}\mathbf{P}\left( S_{k-i}-S_{k}\in A,\,%
S_{k}<L_{k-1},\,S_{k}\leq L_{k+1,n}\right)
\end{eqnarray*}%
and by the Markov property and the duality property of random walks we have
that%
\begin{eqnarray*}
&&\mathbf{P}\left( S_{k-i}-S_{k}\in A,\,S_{k}<L_{k-1},\,%
S_{k}\leq L_{k+1,n}\right) \\
&=&\mathbf{P}\left( S_{k-i}-S_{k}\in A,\,S_{k}<L_{k-1}\right) \mathbf{P%
}\left( S_{k}\leq L_{k+1,n}\right) \\
&=&\mathbf{P}\left( -S_{i}\in A,\,M_{k}\,<0\right) \mathbf{P}%
\left( S_{k}\leq L_{k+1,n}\right) \\
&=&\mathbf{P}\left( \left. -S_{i}\in A\,\right\vert \,M_{k}<0\right) \mathbf{P}\left( M_{k}\,<0\right) \mathbf{P}\left(
S_{k}\leq L_{k+1,n}\right) \\
&=&\mathbf{P}\left( \left. -S_{i}\in A\,\right\vert \,M_{k}<0\right) \mathbf{P}\left( \tau _{n}=k\right) .
\end{eqnarray*}%
Thus,%
\begin{equation}
\mathbf{P}\left( S_{-i,n}^{\prime }\in A,\,\tau _{n}-i\geq 0\right)
=\sum\limits_{k=i}^{n}\mathbf{P}\left( \left. -S_{i}\in A\,\right\vert
\,M_{k}\,<0\right) \mathbf{P}\left( \tau _{n}=k\right) ,
\label{13}
\end{equation}%
therefore, if $\varepsilon \in \left( 0,1\right) $ and $n$ is large enough,
then
\[
\mathbf{P}\left( S_{-i,n}^{\prime }\in A,\,\tau _{n}-i\geq 0\right)
=P_{3}\left( n,\varepsilon \right) +P_{4}\left( n,\varepsilon \right) ,
\]%
where%
\[
P_{3}\left( n,\varepsilon \right) =\sum\limits_{k=i}^{\left\lfloor
\varepsilon n\right\rfloor }\mathbf{P}\left( \left. -S_{i}\in A\,%
\right\vert \,M_{k}\,<0\right) \mathbf{P}\left( \tau
_{n}=k\right) ,
\]%
\[
P_{4}\left( n,\varepsilon \right) =\sum\limits_{k=\left\lfloor \varepsilon
n\right\rfloor +1}^{n}\mathbf{P}\left( \left. -S_{i}\in A\,\right\vert
\,M_{k}\,<0\right) \mathbf{P}\left( \tau _{n}=k\right) .
\]%
It is not difficult to show (see our proof of relation (10)) that

\[
{\lim_{\varepsilon \rightarrow 0}}{%
\limsup_{n\rightarrow \infty }}P_{3}\left( n,\varepsilon \right) =0.
\]

Due to (4) the probability $\mathbf{P}\left( \left. -S_{-i}\in A\,%
\right\vert \,M_{k}\,<0\right) $ tends, as $n\rightarrow \infty $%
, to $\mathbf{P}^{\ast }\left( -S_{-i}^{-}\in A\,\right) =\mathbf{P}%
^{\ast }\left( S_{-i}^{\ast }\in A\,\right) $ uniformly over $%
\left\lfloor \varepsilon n\right\rfloor <k\leq n$. Therefore

\[
{\lim_{\varepsilon \rightarrow 0} }{%
\lim_{n\rightarrow \infty }}P_{4}\left( n,\varepsilon \right)
=\mathbf{P}^{\ast }\left( S_{-i}^{\ast }\in A\,\right) .
\]%
As result, we obtain that

\[
{\lim_{n\rightarrow \infty } }\mathbf{P}\left( S_{-i,n}^{\prime }\in
A,\,\tau _{n}-i\geq 0\right) =\mathbf{P}^{\ast }\left( S_{-i}^{\ast
}\in A\,\right)
\]%
proving (12) for $i\in \mathbf{Z\backslash N}_{0}$. Thus, convergence of
one-dimensional distributions in (6) is established.

The lemma is proved.\medskip

\textbf{Remark 1}. It is not difficult to verify (see \cite{Afanasyev199345}%
, Lemma 1) that relation (6) admits the following generalization: for any $%
a\leq 0$ and $b>0$, as $n\rightarrow \infty $,%
\[
\left\{ S_{i,n}^{\prime },\,i\in \mathbf{Z}\,\left\vert \,%
\frac{L_{n}}{C_{n}}\leq a,\,\frac{S_{n}-L_{n}}{C_{n}}\leq b\right.
\right\} \stackrel{D}{\to}\left\{ S_{i}^{\ast },\,i\in \mathbf{%
Z}\right\} .
\]%
\medskip

Recall that $\left( \Omega ,\mathcal{F},\mathbf{P}\right) $ is the
underlying probability space. Set%
\[
I_{n}^{\left( 2\right) }:=\left\{ \left( i,j\right) :i,j\in \left\{ 0,\ldots
,n\right\} \hbox{ and }i\leq j\right\} .
\]%
Let $\mathcal{F}_{n}$, $n\in \mathbf{N}$, denote the $\sigma $-algebra
generated by the segment of the random environment $Q_{1},\ldots ,Q_{n}$ and
the random variables $Z_{i,j}$ for $\left( i,j\right) \in I_{n}^{\left(
2\right) }$. We now introduce a probability measure $\mathbf{P}^{+}$ on the $%
\sigma $-algebra $\mathcal{F}_{\infty }:=\sigma \left( \cup _{n=1}^{\infty }%
\mathcal{F}_{n}\right) $, defined for each $n\in \mathbf{N}_{0}$ and each $%
\mathcal{F}_{n}$-measurable nonnegative random variable $\beta $ by the
formula%
\begin{equation}
\mathbf{E}^{+}\beta =\mathbf{E}\left( \beta v\left( S_{n}\right) ;\,%
L_{n}\geq 0\right) .  \label{14}
\end{equation}%
This may require a change of the underlying probability space (see \cite%
{Afanasyev2005645} for more details). Similarly, we also introduce a
probability measure $\mathbf{P}^{-}$ on the $\sigma $-algebra $\mathcal{F}%
_{\infty }$, defined for each $n\in \mathbf{N}_{0}$ and each $\mathcal{F}%
_{n} $-measurable nonnegative random variable $\beta $ by the formula

\begin{equation}
\mathbf{E}^{-}\beta =\mathbf{E}\left( \beta u\left( -S_{n}\right) ;\,%
M_{n}<0\right) .  \label{15}
\end{equation}%
Recall that the functions $v\left( \cdot \right) $ and $u\left( \cdot
\right) $ in formulae (14) and (15) are defined after relation (2). Thus,
three measures $\mathbf{P},\mathbf{P}^{+},\mathbf{P}^{-}$ are defined on one
and the same measurable space $\left( \Omega ,\mathcal{F}_{\infty }\right) $%
. To explicitly indicate the measure on $\left( \Omega ,\mathcal{F}_{\infty
}\right) $ according to which we consider this or those random elements we
use the measure symbol as a lower index.

For instance, it is shown in Lemma 2.5 from \cite{Afanasyev2005645} that
\begin{equation}
\left\{ \left( Q_{i}^{+},S_{i}^{+},\mu _{i}^{+}\right) ,\,i\in \mathbf{%
N}\right\} \stackrel{D}{=}\left\{ \left( Q_{i},S_{i},\mu _{i}\right) ,\,%
i\in \mathbf{N}\right\} _{\mathbf{P}^{+}}  \label{16}
\end{equation}%
(the lower index $\mathbf{P}^{+}$ for a random sequence shows here that the
measure $\mathbf{P}^{+}$ is used on the space $\left( \Omega ,\mathcal{F}%
_{\infty }\right) $). Similarly,%
\begin{equation}
\left\{ \left( Q_{i}^{-},S_{i}^{-},\mu _{i}^{-}\right) ,\,i\in \mathbf{%
N}\right\} \stackrel{D}{=}\left\{ \left( Q_{i},S_{i},\mu _{i}\right) ,\,%
i\in \mathbf{N}\right\} _{\mathbf{P}^{-}}.  \label{17}
\end{equation}%
Due to (16), (17) and our assumption about the independence of the left-hand
sides of these relations, the product of probability spaces $\left( \Omega ,%
\mathcal{F}_{\infty },\mathbf{P}^{+}\right) $ and $\left( \Omega ,\mathcal{F}%
_{\infty },\mathbf{P}^{-}\right) $ may be considered as a probability space $%
\left( \Omega ^{\ast },\mathcal{F}^{\ast },\mathbf{P}^{\ast }\right) $ and,
consequently, the direct product of the measures $\mathbf{P}^{+}$ and $%
\mathbf{P}^{-}$ may be treated as the measure $\mathbf{P}^{\ast }$.\medskip

\textbf{Remark 2}. If a random element $\xi $ is given on the space $\left(
\Omega ,\mathcal{F}_{\infty },\mathbf{P}^{+}\right) $ we can define the
random element $\xi ^{+}$, specified on the product of the spaces $\left(
\Omega ,\mathcal{F}_{\infty },\mathbf{P}^{+}\right) $ and $\left( \Omega ,%
\mathcal{F}_{\infty },\mathbf{P}^{-}\right) $ by means of the formula $\xi
^{+}\left( \omega _{1},\omega _{2}\right) =\xi \left( \omega _{1}\right) $
for $\left( \omega _{1},\omega _{2}\right) \in \Omega \times \Omega $. It is
clear that $\mathbf{P}^{\ast }\left( \xi ^{+}\in A\right) =\mathbf{P}%
^{+}\left( \xi \in A\right) $ for an arbitrary one-dimensional Borel set $A$%
. Similarly, if a random element $\xi $ is given on the space $\left( \Omega
,\mathcal{F}_{\infty },\mathbf{P}^{-}\right) $ we can define the random
element $\xi ^{-}$, specified on the product of the spaces $\left( \Omega ,%
\mathcal{F}_{\infty },\mathbf{P}^{+}\right) $ and $\left( \Omega ,\mathcal{F}%
_{\infty },\mathbf{P}^{-}\right) $ by means of the formula $\xi ^{-}\left(
\omega _{1},\omega _{2}\right) =\xi \left( \omega _{2}\right) $ for $\left(
\omega _{1},\omega _{2}\right) \in \Omega \times \Omega $, and $\mathbf{P}%
^{\ast }\left( \xi ^{-}\in A\right) =\mathbf{P}^{-}\left( \xi \in A\right) $
for an arbitrary one-dimensional Borel set $A$. In accordance with the
agreement we can consider the random elements standing in the left-hand
sides of formulae (16) and (17) as generated by the random elements $\left\{
\left( Q_{i},S_{i},\mu _{i}\right) ,\,i\in \mathbf{N}\right\} _{%
\mathbf{P}^{+}}$ and $\left\{ \left( Q_{i},S_{i},\mu _{i}\right) ,\,%
i\in \mathbf{N}\right\} _{\mathbf{P}^{-}}$ respectively.\medskip

\textbf{Lemma 2}.\ \textit{If the Spitzer-Doney condition} (2) \textit{is
satisfied}, \textit{then}, \textit{as} $n\rightarrow \infty $,%
\begin{equation}
\left\{ \left. a_{i,n}Z_{i,n},\,i\in \mathbf{N}_{0}\,\right\vert
\,L_{n}\,\geq 0\right\} \stackrel{D}{\to}\left\{ \zeta _{i}^{\ast
},\,i\in \mathbf{N}_{0}\right\} ,  \label{18}
\end{equation}%
\textit{where} $\left\{ \zeta _{i}^{\ast },\,i\in \mathbf{N}%
_{0}\right\} $ \textit{is the random sequence defined by relation} (5)%
\textit{.\medskip }

\textit{Proof}. By virtue of the first part of Lemma 2.5 from \cite%
{Afanasyev2005645} for $k\in \mathbf{N}$, as $n\rightarrow \infty $,%
\[
\left\{ \left. \left( a_{i,j},Z_{i,j}\right) ,\,\left( i,j\right) \in
I_{k}^{\left( 2\right) }\,\right\vert \,L_{n}\,\geq
0\right\} \stackrel{D}{\to}\left\{ \left( a_{i,j},Z_{i,j}\right) ,%
\,\left( i,j\right) \in I_{k}^{\left( 2\right) }\right\} _{\mathbf{P}%
^{+}}.
\]%
Note (see \cite{Afanasyev2005645}, Section 3) that in view of (14) for a
fixed $i\in \mathbf{N}_{0}$ the random sequence $\left\{ \eta _{i};\,%
Z_{i,j},\,j\in \mathbf{N}_{i+1}\right\} _{\mathbf{P}^{+}}$ given on
the probability space $\left( \Omega ,\mathcal{F}_{\infty },\mathbf{P}%
^{+}\right) $ is a branching process in the random environment $\left\{
G_{i+1};\,F_{n},\,n\in \mathbf{N}_{i+1}\right\} _{\mathbf{P}%
^{+}} $. Hence, if the random environment is fixed, the sequence $\left\{
\eta _{i};\,a_{i,j}Z_{i,j},\,j\in \mathbf{N}_{i+1}\right\} _{%
\mathbf{P}^{+}}$ is a non-negative martingale. Because of this (without
assuming that the random environment is fixed) there is $\mathbf{P}^{+}$%
-a.s. the finite limit%
\[
\lim_{n\rightarrow \infty }a_{i,n}Z_{i,n}=:\zeta _{i}.
\]%
It means, in view of the second part of Lemma 2.5 from \cite%
{Afanasyev2005645}, that, as $n\rightarrow \infty $,

\[
\left\{ \left. a_{i,n}Z_{i,n},\,i\in \mathbf{N}_{0}\,\right\vert
\,L_{n}\,\geq 0\right\} \stackrel{D}{\to}\left\{ \zeta _{i},\,i\in
\mathbf{N}_{0}\right\} _{\mathbf{P}^{+}}.
\]%
To prove relation (18), it remains to note that in view of Remark 2%
\[
\left\{ \zeta _{i},\,i\in \mathbf{N}_{0}\right\} _{\mathbf{P}^{+}}%
\stackrel{D}{=}\left\{ \zeta _{i}^{\ast },\,i\in
\mathbf{N}_{0}\right\} .
\]

The lemma is proved.\medskip

Set for $n\in \mathbf{N}$%
\[
Z_{i,n}^{\prime }=\left\{
\begin{array}{c}
Z_{\tau _{n}+i,n},\qquad i\in \mathbf{N}_{\left( -\tau _{n}\right) },
\\
0,\qquad i\in \mathbf{Z\setminus N}_{\left( -\tau _{n}\right) },%
\end{array}%
\right.
\]%
\[
a_{i,n}^{\prime }=\left\{
\begin{array}{c}
a_{n}/a_{\tau _{n}+i},\qquad i\in \mathbf{N}_{\left( -\tau
_{n}\right) }, \\
1,\qquad i\in \mathbf{Z\setminus N}_{\left( -\tau _{n}\right) }.%
\end{array}%
\right.
\]%
\medskip

\textbf{Lemma 3}.\ \textit{If the Spitzer-Doney condition} (2) \textit{is
satisfied}, \textit{then}, \textit{as} $n\rightarrow \infty $,%
\begin{equation}
\left\{ a_{i,n}^{\prime }Z_{i,n}^{\prime },\,i\in \mathbf{Z}\right\}
\stackrel{D}{\to}\left\{ \zeta _{i}^{\ast },\,i\in \mathbf{Z}%
\right\} ,  \label{19}
\end{equation}%
\textit{where} $\left\{ \zeta _{i}^{\ast },\,i\in \mathbf{N}%
_{0}\right\} $ \textit{is the random sequence defined by relation}
(5).\medskip

\textit{Proof}. We demonstrate for simplicity only convergence of
one-dimensional distributions. Fix $i\in \mathbf{N}_{0}$. Let $A$ be an
arbitrary one-dimensional $\zeta _{i}^{\ast }$-continuous (relative to the
measure $\mathbf{P}^{\ast }$) Borel set. Then for $n\geq i$%
\begin{eqnarray*}
&&\mathbf{P}\left( a_{i,n}^{\prime }Z_{i,n}^{\prime }\in A,\,\tau
_{n}+i\leq n\right) \\
&=&\sum\limits_{k=0}^{n-i}\mathbf{P}\left( a_{i,n}^{\prime }Z_{i,n}^{\prime
}\in A,\,\tau _{n}=k\right) \\
&=&\sum\limits_{k=0}^{n-i}\mathbf{P}\left( \frac{a_{n}}{a_{k+i}}Z_{k+i,n}\in
A,\,S_{k}<L_{k-1},\,S_{k}\leq L_{k+1,n}\right)
\end{eqnarray*}%
and

\begin{eqnarray*}
&&\mathbf{P}\left( \frac{a_{n}}{a_{k+i}}Z_{k+i,n}\in A,\,S_{k}<L_{k-1},%
\,S_{k}\leq L_{k+1,n}\right) \\
&=&\mathbf{P}\left( S_{k}<L_{k-1}\right) \mathbf{P}\left( \frac{a_{n}}{%
a_{k+i}}Z_{k+i,n}\in A,\,S_{k}\leq L_{k+1,n}\right) \\
&=&\mathbf{P}\left( S_{k}<L_{k-1}\right) \mathbf{P}\left( \frac{a_{n-k}}{%
a_{i}}Z_{i,n-k}\in A,\,L_{n-k}\,\geq 0\right) \\
&=&\mathbf{P}\left( \left. a_{i,n-k}Z_{i,n-k}\in A\,\right\vert\, L_{n-k}\,\geq 0\right) \mathbf{P}\left( S_{k}<L_{k-1}\right) \mathbf{P%
}\left( L_{n-k}\,\geq 0\right) \\
&=&\mathbf{P}\left( \left. a_{i,n-k}Z_{i,n-k}\in A\,\right\vert \,L_{n-k}\,\geq 0\right) \mathbf{P}\left( \tau _{n}=k\right) .
\end{eqnarray*}%
Thus,%
\begin{eqnarray*}
&&\mathbf{P}\left( a_{i,n}^{\prime }Z_{i,n}^{\prime }\in A,\,\tau
_{n}+i\leq n\right) \\
&=&\sum\limits_{k=0}^{n-i}\mathbf{P}\left( \left. a_{i,n-k}Z_{i,n-k}\in A%
\,\right\vert \,L_{n-k}\,\geq 0\right) \mathbf{P}\left(
\tau _{n}=k\right) .
\end{eqnarray*}%
Therefore, if $\varepsilon \in \left( 0,1\right) $ and $n$ is large enough,
then
\[
\mathbf{P}\left( a_{i,n}^{\prime }Z_{i,n}^{\prime }\in A,\,\tau
_{n}+i\leq n\right) =P_{1}\left( n,\varepsilon \right) +P_{2}\left(
n,\varepsilon \right) ,
\]%
where%
\[
P_{1}\left( n,\varepsilon \right) =\sum\limits_{k=0}^{\left\lfloor \left(
1-\varepsilon \right) n\right\rfloor }\mathbf{P}\left( \left.
a_{i,n-k}Z_{i,n-k}\in A\,\right\vert \,L_{n-k}\,\geq
0\right) \mathbf{P}\left( \tau _{n}=k\right) ,
\]%
\[
P_{2}\left( n,\varepsilon \right) =\sum\limits_{k=\left\lfloor \left(
1-\varepsilon \right) n\right\rfloor +1}^{n-i}\mathbf{P}\left( \left.
a_{i,n-k}Z_{i,n-k}\in A\,\right\vert \,L_{n-k}\,\geq
0\right) \mathbf{P}\left( \tau _{n}=k\right) .
\]%
It is easy to show (see the proof of relation (10)) that%
\[
{\lim_{\varepsilon \rightarrow 0} }{%
\limsup_{n\rightarrow \infty } }P_{2}\left( n,\varepsilon \right)
=0.
\]%
By Lemma 2 the probability $\mathbf{P}\left( \left. a_{i,n-k}Z_{i,n-k}\in A%
\,\right\vert \,L_{n-k}\,\geq 0\right) $ tends, as $%
n\rightarrow \infty $, to $\mathbf{P}^{\ast }\left( \zeta _{i}^{\ast }\in A%
\,\right) $ uniformly over $0\leq k\leq \left\lfloor \left(
1-\varepsilon \right) n\right\rfloor $. Therefore%
\[
{\lim_{\varepsilon \rightarrow 0} }{%
\lim_{n\rightarrow \infty } }P_{1}\left( n,\varepsilon \right)
=\mathbf{P}^{\ast }\left( \zeta _{i}^{\ast }\in A\,\right) .
\]%
As result, we obtain that%
\[
{\lim_{n\rightarrow \infty } }\mathbf{P}\left( a_{i,n}^{\prime
}Z_{i,n}^{\prime }\in A,\,\tau _{n}+i\leq n\right) =\mathbf{P}^{\ast
}\left( \zeta _{i}^{\ast }\in A\,\right) .
\]%
This justifies the one-dimensional convergence in (19) for $i\in \mathbf{N}%
_{0}$.

Now fix $i\in \mathbf{N}$. Then for $x\geq 0$ and $n\geq i$%
\begin{eqnarray}
&&\mathbf{P}\left( a_{-i,n}^{\prime }Z_{-i,n}^{\prime }\leq x,\,\tau
_{n}-i\geq 0\right)  \nonumber \\
&=&\sum\limits_{k=i}^{n}\mathbf{P}\left( a_{-i,n}^{\prime }Z_{-i,n}^{\prime
}\leq x,\,\tau _{n}=k\right) \nonumber\\
&=&\sum\limits_{k=i}^{n}\mathbf{P}\left( \frac{a_{n}}{a_{k-i}}Z_{k-i,n}\leq
x,\,S_{k}<L_{k-1},\,S_{k}\leq L_{k+1,n}\right) .  \label{20}
\end{eqnarray}

Note that the random sequence $\left\{ \left( Z_{k-i,n},a_{k-i,n}\right) ,%
\,n\in \mathbf{N}_{k-i}\right\} $ is Markovian. Denote by $%
Z_{k,n}\left( l\right) $ the number of particles of $n$th generation
being descendants of $l$ particles of $k$th generation. Since
$Z_{k-i,n}\stackrel{D}{=}%
Z_{k,n}\left( l\right) $ given $Z_{k-i,k}=l$, it follows that%
\begin{eqnarray*}
&&\mathbf{P}\left( \frac{a_{n}}{a_{k-i}}Z_{k-i,n}\leq x,\,%
S_{k}<L_{k-1},\,S_{k}\leq L_{k+1,n}\right) \\
&=&\mathbf{E}\left( U\left( Z_{k-i,k},a_{k-i,k}\right) ;\,%
S_{k}<L_{k-1}\right) ,
\end{eqnarray*}%
where%
\[
U\left( l,y\right) =\mathbf{P}\left( a_{k,n}Z_{k,n}\left( l\right) \leq
\frac{x}{y},\,S_{k}\leq L_{k+1,n}\right) .
\]%
Clearly,%
\[
U\left( l,y\right) =\mathbf{P}\left( a_{0,n-k}Z_{0,n-k}\left( l\right) \leq
\frac{x}{y},\,L_{n-k}\,\geq 0\right) .
\]%
As result, we obtain that%
\begin{eqnarray}
&&\mathbf{P}\left( \frac{a_{n}}{a_{k-i}}Z_{k-i,n}\leq x,\,%
S_{k}<L_{k-1},\,S_{k}\leq L_{k+1,n}\right) \nonumber\\
&=&\mathbf{E}\left( \left. H_{n-k}\left( Z_{k-i,k},\frac{x}{a_{k-i,k}}%
\right) \,\right\vert \,S_{k}<L_{k-1}\right) \mathbf{P}\left(
S_{k}<L_{k-1}\right) \mathbf{P}\left( L_{n-k}\,\geq 0\right) \nonumber\\
&=&\mathbf{E}\left( \left. H_{n-k}\left( Z_{k-i,k},x/a_{k-i,k}\right)\, \right\vert \,S_{k}<L_{k-1}\right) \mathbf{P}\left( \tau
_{n}=k\right) ,  \label{21}
\end{eqnarray}%
where%
\[
H_{n}\left( l,x\right) =\mathbf{P}\left( \left. a_{0,n}Z_{0,n}\left(
l\right) \leq x\,\right\vert \,L_{n}\,\geq 0\right)
\]%
for $l\in \mathbf{N}_{0}$ and $x\geq 0$.

Set $Q_{k,l}=\left( Q_{k},\ldots ,Q_{l}\right) $ for $k,l\in \mathbf{N}$. In
the sequel, we will to explicitly include a random environment in the
notation. For example, we will write $Z_{k-i,k}\left\langle
Q_{k-i+1,k}\right\rangle $ instead of $Z_{k-i,k}$. Set%
\[
Q_{k}=\widetilde{Q}_{1},\ldots ,Q_{k-i+1}=\widetilde{Q}_{i},\ldots ,Q_{1}=%
\widetilde{Q}_{k}
\]%
and consider a branching process with immigration in the random environment $%
\widetilde{Q}_{1},\ldots ,\widetilde{Q}_{k}$. Then%
\[
\begin{array}{c}
\mathbf{E}\left( \left. H_{n-k}\left( Z_{k-i,k}\left\langle
Q_{k-i+1,k}\right\rangle ,x/a_{k-i,k}\right) \,\right\vert \,%
S_{k}<L_{k-1}\right) \qquad\qquad \\
=\mathbf{E}\left( \left. H_{n-k}\left( Z_{0,i}\left\langle \widetilde{Q}%
_{i,1}\right\rangle ,x/\widetilde{a}_{0,i}\right) \,\right\vert\, \widetilde{M}_{k}\,<0\right) ,%
\end{array}%
\]%
where the symbols $\widetilde{a}_{0,i},\widetilde{M}_{k},\widetilde{Q}_{i,1}$
have the same meaning for the random environment $\widetilde{Q}_{1},,\ldots ,%
\widetilde{Q}_{k}$ as the symbols $a_{0,i},M_{k},Q_{i,1}$ mean for the
random environment $Q_{1},\ldots ,Q_{k}$. Further, the random environments $%
\widetilde{Q}_{1},,\ldots ,\widetilde{Q}_{k}$ and $Q_{1},,\ldots ,Q_{k}$ are
identically distributed. Therefore%
\begin{eqnarray*}
&&\mathbf{E}\left( \left. H_{n-k}\left( Z_{0,i}\left\langle \widetilde{Q}%
_{i,1}\right\rangle ,x/\widetilde{a}_{0,i}\right) \,\right\vert\, \widetilde{M}_{k}\,<0\right) \\
&=&\mathbf{E}\left( \left. H_{n-k}\left( Z_{0,i}\left\langle
Q_{i,1}\right\rangle ,x/a_{0,i}\right) \,\right\vert \,M_{k}%
\,<0\right) .
\end{eqnarray*}%
As result, we obtain that%
\begin{eqnarray}
&&\mathbf{E}\left( \left. H_{n-k}\left( Z_{k-i,k}\left\langle
Q_{k-i+1,k}\right\rangle ,x/a_{k-i,k}\right) \,\right\vert \,%
S_{k}<L_{k-1}\right) \nonumber\\
&=&\mathbf{E}\left( \left. H_{n-k}\left( Z_{0,i}\left\langle
Q_{i,1}\right\rangle ,x/a_{0,i}\right) \,\right\vert \,M_{k}%
\,<0\right) .  \label{22}
\end{eqnarray}

Set $\psi _{i}=Z_{0,i}\left\langle Q_{i,1}\right\rangle $. We have from
(20)-(22) that%
\begin{eqnarray*}
&&\mathbf{P}\left( a_{-i,n}^{\prime }Z_{-i,n}^{\prime }\leq x,\,\tau
_{n}-i\geq 0\right) \\
&=&\sum\limits_{k=i}^{n}\mathbf{E}\left( \left. H_{n-k}\left( \psi
_{i},x/a_{0,i}\right) \,\right\vert \,M_{k}\,<0\right)
\mathbf{P}\left( \tau _{n}=k\right) .
\end{eqnarray*}%
Therefore, if $\varepsilon \in \left( 0,1\right) $ and $n$ is large enough,
then%
\begin{equation}
\mathbf{P}\left( a_{-i,n}^{\prime }Z_{-i,n}^{\prime }\leq x,\,\tau
_{n}-i\geq 0\right) =P_{3}\left( n,\varepsilon \right) +P_{4}\left(
n,\varepsilon \right) +P_{5}\left( n,\varepsilon \right) ,  \label{23}
\end{equation}%
where%
\[
P_{3}\left( n,\varepsilon \right) =\sum\limits_{k=i}^{\left\lfloor
\varepsilon n\right\rfloor }\mathbf{E}\left( \left. H_{n-k}\left( \psi
_{i},x/a_{0,i}\right) \,\right\vert \,M_{k}\,<0\right)
\mathbf{P}\left( \tau _{n}=k\right) ,
\]%
\[
P_{4}\left( n,\varepsilon \right) =\sum\limits_{k=\left\lfloor \left(
1-\varepsilon \right) n\right\rfloor +1}^{n}\mathbf{E}\left( \left.
H_{n-k}\left( \psi _{i},x/a_{0,i}\right) \,\right\vert \,M_{k}%
\,<0\right) \mathbf{P}\left( \tau _{n}=k\right) ,
\]%
\[
P_{5}\left( n,\varepsilon \right) =\sum\limits_{k=\left\lfloor \varepsilon
n\right\rfloor +1}^{\left\lfloor \left( 1-\varepsilon \right) n\right\rfloor
}\mathbf{E}\left( \left. H_{n-k}\left( \psi _{i},x/a_{0,i}\right) \,%
\right\vert \,M_{k}\,<0\right) \mathbf{P}\left( \tau
_{n}=k\right) .
\]

Similar to relation (10) we conclude that%
\begin{equation}
{\lim_{\varepsilon \rightarrow 0} }{%
\limsup_{n\rightarrow \infty } }P_{3}\left( n,\varepsilon \right)
=0, \label{24}
\end{equation}%
\begin{equation}
{\lim_{\varepsilon \rightarrow 0} }{%
\limsup_{n\rightarrow \infty } }P_{4}\left( n,\varepsilon \right)
=0. \label{25}
\end{equation}

Let $l\in \mathbf{N}_{0}$ be fixed. It is not difficult to demonstrate that%
\begin{equation}
\lim_{n\rightarrow \infty }a_{0,n}Z_{0,n}\left( l\right) =:\zeta _{0}\left(
l\right)  \label{26}
\end{equation}%
exists a.s. on the probability space $\left( \Omega ,\mathcal{F}_{\infty },%
\mathbf{P}^{+}\right) $. By the arguments to those used in Lemma 2 one can
show, as $n\rightarrow \infty $,%
\begin{equation}
\left\{ \left. a_{0,n}Z_{0,n}\left( l\right) ,\,i\in \mathbf{N}_{0}%
\,\right\vert \,L_{n}\,\geq 0\right\} \stackrel{D}{\to}%
\left\{ \zeta _{0}\left( l\right) ,\,i\in \mathbf{N}%
_{0}\right\} _{\mathbf{P}^{+}}.  \label{27}
\end{equation}%
For $x\geq 0$ set%
\[
H\left( l,x\right) =\mathbf{P}^{+}\left( \zeta _{0}\left( l\right) \leq
x\right) .
\]%
It follows from (27) that%
\begin{equation}
\lim_{n\rightarrow \infty }H_{n}\left( l,x\right) =H\left( l,x\right)
\label{28}
\end{equation}%
if $x\geq 0$ belongs to the set of continuity points of $H\left( l,\cdot
\right) $ (with respect to the second argument). By Lemma 2.5 in \cite%
{Afanasyev2005645}%
\[
\left\{ \left. Z_{0,i}\left\langle Q_{i,1}\right\rangle ,a_{0,i}\,%
\right\vert \,M_{n}\,<0\right\} \stackrel{D}{\to}\left(
Z_{0,i}\left\langle Q_{i,1}\right\rangle ,a_{0,i}\right)
_{\mathbf{P}^{-}}
\]%
as $n\rightarrow \infty $. Therefore%
\begin{equation}
\left\{ \left. \psi _{i},a_{0,i}\,\right\vert \,M_{n}\,%
<0\right\} \stackrel{D}{\to}\left( Z_{-i,0}^{\ast },a_{-i,0}^{\ast
}\right) . \label{29}
\end{equation}

We show that, for fixed $l\in \mathbf{N}_{0}$ and $K>0$%
\begin{eqnarray}
&&\lim_{n\rightarrow \infty }\mathbf{E}\left( \left. H_{n-k}\left(
l,x/a_{0,i}\right) I_{\left\{ \psi _{i}=l,\,x/a_{0,i}\leq K\right\} }%
\,\right\vert \,M_{k}\,<0\right) \nonumber\\
&=&\mathbf{E}^{\ast }\left( H\left( l,x/a_{-i,0}^{\ast }\right) ;\,%
Z_{-i,0}^{\ast }=l,\,x/a_{-i,0}^{\ast }\leq K\right)  \label{30}
\end{eqnarray}%
uniformly over $\left\lfloor \varepsilon n\right\rfloor <k\leq \left\lfloor
\left( 1-\varepsilon \right) n\right\rfloor $ (here $I_{A}$ is the indicator
of the event $A$). Let $0=x_{0}<x_{1}<\ldots <x_{m}=K$ for some $m\in
\mathbf{N}$. The monotonicity of the function $H\left( l,\cdot \right) $
with respect to the second argument gives%
\begin{eqnarray}
&&\mathbf{E}\left( \left. H_{n-k}\left( l,x/a_{0,i}\right) I_{\left\{ \psi
_{i}=l,\,x/a_{0,i}\leq K\right\} }\,\right\vert \,M_{k}%
\,<0\right) \nonumber\\
&=&\sum\limits_{j=1}^{m}\mathbf{E}\left( \left. H_{n-k}\left(
l,x/a_{0,i}\right) I_{\left\{ \psi _{i}=l,\,x_{j-1}<x/a_{0,i}\leq
x_{j}\right\} }\,\right\vert \,M_{k}\,<0\right) \nonumber\\
&\leq &\sum\limits_{j=1}^{m}H_{n-k}\left( l,x_{j}\right) \mathbf{P}\left(
\left. \psi _{i}=l,\,x_{j-1}<x/a_{0,i}\leq x_{j}\,\right\vert
\,M_{k}\,<0\right) .  \label{31}
\end{eqnarray}%
In view of (28) and (29) the right-hand side of (31) converges, as $%
n\rightarrow \infty $, to%
\[
\sum\limits_{j=1}^{m}H\left( l,x_{j}\right) \mathbf{P}^{\ast }\left(
Z_{-i,0}^{\ast }=l,\,x_{j-1}<x/a_{-i,0}^{\ast }\leq x_{j}\right)
\]%
uniformly over $\left\lfloor \varepsilon n\right\rfloor <k\leq \left\lfloor
\left( 1-\varepsilon \right) n\right\rfloor $, if the selected $x_{1},\ldots
,x_{m}$ are simultaneously the continuity points of $H\left( l,\cdot \right)
$ with respect to the second argument and of $\mathbf{P}^{\ast }\left(
Z_{-i,0}^{\ast }=l,\,x/a_{-i,0}^{\ast }\leq y\right) $ with respect to
$y$. Thus, if $\delta >0$ and $n$ is large enough, the following inequality
holds%
\begin{eqnarray}
&&\mathbf{E}\left( \left. H_{n-k}\left( l,x/a_{0,i}\right) I_{\left\{ \psi
_{i}=l,\,x/a_{0,i}\leq K\right\} }\,\right\vert \,M_{k}%
\,<0\right) \nonumber\\
&\leq &\sum\limits_{j=1}^{m}H\left( l,x_{j}\right) \mathbf{P}^{\ast }\left(
Z_{-i,0}^{\ast }=l,\,x_{j-1}<\left( a_{-i,0}^{\ast }\right) ^{-1}x\leq
x_{j}\right) +\delta  \label{32}
\end{eqnarray}%
for $\left\lfloor \varepsilon n\right\rfloor <k\leq \left\lfloor \left(
1-\varepsilon \right) n\right\rfloor $. Similarly, if $\delta >0$ and $n$ is
large enough, then%
\begin{eqnarray}
&&\mathbf{E}\left( \left. H_{n-k}\left( l,x/a_{0,i}\right) I_{\left\{ \psi
_{i}=l,\,x/a_{0,i}\leq K\right\} }\,\right\vert \,M_{k}%
\,<0\right) \nonumber\\
&\geq &\sum\limits_{j=1}^{m}H\left( l,x_{j-1}\right) \mathbf{P}^{\ast
}\left( Z_{-i,0}^{\ast }=l,\,x_{j-1}<x/a_{-i,0}^{\ast }\leq
x_{j}\right) -\delta  \label{33}
\end{eqnarray}%
for $\left\lfloor \varepsilon n\right\rfloor <k\leq \left\lfloor \left(
1-\varepsilon \right) n\right\rfloor $. Since $0\leq H\left( l,x\right) \leq
1$ for $x\geq 0$, the sums in the right-hand sides of (32) and (33)
converge, as $\max_{1\leq j\leq m}\left( x_{j}-x_{j-1}\right) \rightarrow 0$%
, to (see \cite{Shiryaev2018}, Chapter 2, \S\ 6, Section 11)%
\[
\mathbf{E}^{\ast }\left[ H\left( l,x/a_{-i,0}^{\ast }\right) ;\,%
Z_{-i,0}^{\ast }=l,\,x/a_{-i,0}^{\ast }\leq K\right] .
\]%
Hence, if $\delta >0$ and $n$ is large enough, then%
\begin{eqnarray*}
&&\mathbf{E}^{\ast }\left( H\left( l,x/a_{-i,0}^{\ast }\right) ;\,%
Z_{-i,0}^{\ast }=l,\,x/a_{-i,0}^{\ast }\leq K\right) -\delta \\
&\leq &\mathbf{E}\left( \left. H_{n-k}\left( l,x/a_{0,i}\right) I_{\left\{
\psi _{i}=l,\,x/a_{0,i}\leq K\right\} }\,\right\vert \,%
M_{k}\,<0\right) \\
&\leq &\mathbf{E}^{\ast }\left( H\left( l,x/a_{-i,0}^{\ast }\right) ;\,%
Z_{-i,0}^{\ast }=l,\,x/a_{-i,0}^{\ast }\leq K\right) +\delta
\end{eqnarray*}%
for $\left\lfloor \varepsilon n\right\rfloor <k\leq \left\lfloor \left(
1-\varepsilon \right) n\right\rfloor $. Since $\delta >0$ is arbitrary, we
obtain the required relation (30).

Now we show that%
\begin{equation}
\lim_{n\rightarrow \infty }\mathbf{E}\left( \left. H_{n-k}\left( \psi
_{i},x/a_{0,i}\right) \,\right\vert \,M_{k}\,<0\right) =%
\mathbf{E}^{\ast }H\left( Z_{-i,0}^{\ast },x/a_{-i,0}^{\ast }\right)
\label{34}
\end{equation}%
uniformly over $\left\lfloor \varepsilon n\right\rfloor <k\leq \left\lfloor
\left( 1-\varepsilon \right) n\right\rfloor $. For $N\in \mathbf{N}$ and $K>0
$ we write%
\begin{equation}
\mathbf{E}\left( \left. H_{n-k}\left( \psi _{i},x/a_{0,i}\right) \,%
\right\vert \,M_{k}\,<0\right) =E_{1}\left( k,n,N,K\right)
+E_{2}\left( k,n,N,K\right) ,  \label{35}
\end{equation}%
where%
\[
E_{1}\left( k,n,N,K\right) =\mathbf{E}\left( \left. H_{n-k}\left( \psi
_{i},x/a_{0,i}\right) I_{\left\{ \psi _{i}\leq N,\,x/a_{0,i}\leq
K\right\} }\,\right\vert \,M_{k}\,<0\right) ,
\]%
\[
E_{2}\left( k,n,N,K\right) =\mathbf{E}\left( \left. H_{n-k}\left( \psi
_{i},x/a_{0,i}\right) I_{\left\{ \psi _{i}>N\right\} \cup \left\{
x/a_{0,i}>K\right\} }\,\right\vert \,M_{k}\,<0\right) .
\]%
Since%
\[
E_{2}\left( k,n,N,K\right) \leq \mathbf{P}\left( \left. \left\{ \psi
_{i}>N\right\} \cup \left\{ x/a_{0,i}>K\right\} \,\right\vert \,%
M_{k}\,<0\right) ,
\]%
it follows by (29) that%
\begin{equation}
{\lim_{K\rightarrow \infty } }{\lim_{N\rightarrow \infty } }%
{\limsup_{n\rightarrow \infty } }E_{2}\left( k,n,N,K\right) =0
\label{36}
\end{equation}%
uniformly over $\left\lfloor \varepsilon n\right\rfloor <k\leq \left\lfloor
\left( 1-\varepsilon \right) n\right\rfloor $.\ Clearly,%
\[
E_{1}\left( k,n,N,K\right) =\sum\limits_{l=0}^{N}\mathbf{E}\left( \left.
H_{n-k}\left( l,x/a_{0,i}\right) I_{\left\{ \psi _{i}=l,\,%
x/a_{0,i}\leq K\right\} }\,\right\vert \,M_{k}\,<0\right) .
\]%
Hence, using (30) we conclude that%
\begin{equation}
{\lim_{K\rightarrow \infty }}{\lim_{N\rightarrow \infty } }%
\lim_{n\rightarrow \infty }E_{1}\left( k,n,N,K\right)
=\mathbf{E}^{\ast }H\left( Z_{-i,0}^{\ast },x/a_{-i,0}^{\ast
}\right) .  \label{37}
\end{equation}%
uniformly over $\left\lfloor \varepsilon n\right\rfloor <k\leq \left\lfloor
\left( 1-\varepsilon \right) n\right\rfloor $. Combining (35)-(37) we obtain
the desired relation (34).

It follows from (34) that%
\begin{equation}
\lim_{\varepsilon \rightarrow 0}\lim_{n\rightarrow \infty }P_{5}\left(
n,\varepsilon \right) =\mathbf{E}^{\ast }H\left( Z_{-i,0}^{\ast
},x/a_{-i,0}^{\ast }\right) .  \label{38}
\end{equation}%
Now (23)-(25) and (38) imply%
\[
\lim_{n\rightarrow \infty }\mathbf{P}\left( a_{-i,n}^{\prime
}Z_{-i,n}^{\prime }\leq x,\,\tau _{n}-i\geq 0\right) =\mathbf{E}^{\ast
}H\left( Z_{-i,0}^{\ast },x/a_{-i,0}^{\ast }\right) .
\]%
Hence,%
\begin{equation}
\lim_{n\rightarrow \infty }\mathbf{P}\left( a_{-i,n}^{\prime
}Z_{-i,n}^{\prime }\leq x\right) =\mathbf{E}^{\ast }H\left( Z_{-i,0}^{\ast
},x/a_{-i,0}^{\ast }\right) .  \label{39}
\end{equation}

We now analyze a branching process with immigration in the random
environment $\left\{ Q_{k}^{\ast },\,k\in \mathbf{Z}\right\} $. The
random sequence $\left\{ \left( Z_{-i,n}^{\ast },a_{-i,n}^{\ast }\right) ,%
\,n\in \mathbf{N}_{-i}\right\} $ is Markovian. Denote by $%
Z_{k,n}^{\ast }\left( l\right) $ the number of particles in $n$th generation
which are descendants of $l$ particles of the $k$th generation. Note that%
\begin{equation}
\left\{ \left( Z_{0,n}^{\ast }\left( l\right) ,a_{0,n}^{\ast
}\right) ,\, n\in \mathbf{N}_{0}\right\} \stackrel{D}{=}\left\{
\left( Z_{0,n}\left(
l\right) ,a_{0,n}\right) ,\,n\in \mathbf{N}_{0}\right\} _{\mathbf{P}%
^{+}}.  \label{40}
\end{equation}%
Since $Z_{-i,n}^{\ast }\stackrel{D}{=}Z_{0,n}^{\ast }\left( l\right) $ given $%
Z_{-i,0}^{\ast }=l$, it follows that, for any bounded and continuous
function $f:\mathbf{R}\rightarrow \mathbf{R}$ and $n\in \mathbf{N}_{-i}$,%
\begin{equation}
\mathbf{E}^{\ast }f\left( a_{-i,n}^{\ast }Z_{-i,n}^{\ast }\right) =\mathbf{E}%
^{\ast }V_{n}\left( Z_{-i,0}^{\ast },1/a_{-i,0}^{\ast }\right) ,  \label{41}
\end{equation}%
where%
\[
V_{n}\left( l,y\right) =\mathbf{E}^{\ast }f\left( a_{0,n}^{\ast
}Z_{0,n}^{\ast }\left( l\right) /y\right) .
\]%
By (40)%
\begin{equation}
V_{n}\left( l,y\right) =\mathbf{E}^{+}f\left( a_{0,n}Z_{0,n}\left( l\right)
/y\right) .  \label{42}
\end{equation}%
In view of (5), as $n\rightarrow \infty $,
\begin{equation}
a_{-i,n}^{\ast }Z_{-i,n}^{\ast }\stackrel{D}{\to}\zeta _{-i}^{\ast
}, \label{43}
\end{equation}%
and in view of (26)%
\begin{equation}
\left( a_{0,n}Z_{0,n}\left( l\right) \right) _{\mathbf{P}^{+}}\stackrel{D}{\to}%
\left( \zeta _{0}\left( l\right) \right) _{\mathbf{P}^{+}}.
\label{44}
\end{equation}%
Using (42), (44) and applying the dominated convergence theorem we see that%
\begin{equation}
\lim_{n\rightarrow \infty }V_{n}\left( l,y\right) =V\left( l,y\right) ,
\label{45}
\end{equation}%
where%
\[
V\left( l,y\right) =\mathbf{E}^{+}f\left( \zeta _{0}\left( l\right)
/y\right) .
\]%
Applying the dominated convergence theorem again we obtain from (41), (43)
and (45) that
\begin{equation}
\mathbf{E}^{\ast }f\left( \zeta _{-i}^{\ast }\right) =\mathbf{E}^{\ast
}V\left( Z_{-i,0}^{\ast },1/a_{-i,0}^{\ast }\right) .  \label{46}
\end{equation}%
Fix $x\geq 0$. As relation (46) is valid for any bounded and continuous
function $f$, it is valid, even when a function $f$ is the indicator of the
semi-axis $\left( -\infty ,x\right] $. It means that%
\begin{equation}
\mathbf{P}^{\ast }\left( \zeta _{-i}^{\ast }\leq x\right) =\mathbf{E}^{\ast
}H\left( Z_{-i,0}^{\ast },x/a_{-i,0}^{\ast }\right)  \label{47}
\end{equation}%
(we take into account that $V\left( l,y\right) =H\left( l,xy\right) $ for
the specified function $f$).

Equalities (39) and (47) imply the one-dimensional convergence in relation
(19) for $i\in \mathbf{Z\setminus N}_{0}$.

The lemma is proved.\medskip

\textbf{Remark 3}. It is not difficult to verify that (19) admits the
following generalization: for any $a\leq 0$ and $b>0$, as $n\rightarrow
\infty $,%
\[
\left\{ a_{i,n}^{\prime }Z_{i,n}^{\prime },\,i\in \mathbf{Z}\,%
\left\vert \,\frac{L_{n}}{C_{n}}\leq a,\,\frac{S_{n}-L_{n}}{C_{n}%
}\leq b\right. \right\} \stackrel{D}{\to}\left\{ \zeta _{i}^{\ast },%
\,i\in \mathbf{Z}\right\} .
\]%
\medskip

\textbf{Lemma 4}. \textit{If the conditions} \textit{of Theorem 1 are}
\textit{satisfied}, \textit{then} $\mathbf{P}^{\ast }$-\textit{a.s.}%
\[
\Sigma _{1}<+\infty ,\qquad\Sigma _{2}<+\infty .
\]%
\medskip

\textit{Proof. }It is shown in Lemma 2.7 from \cite{Afanasyev2005645} that,
if\ the conditions of Theorem 1 are satisfied, then the series $%
\sum\nolimits_{i=0}^{\infty }\mu _{i+1}e^{-S_{i}}\ $converges $\mathbf{P}%
^{+} $-a.s. Hence, the series $\sum\nolimits_{i=0}^{\infty }\mu
_{i+1}^{+}e^{-S_{i}^{+}}$ converges $\mathbf{P}^{\ast }$-a.s. Similarly we
can prove that the series $\sum\nolimits_{i=1}^{\infty }\mu
_{i}^{-}e^{S_{i}^{-}}$ converges $\mathbf{P}^{\ast }$-a.s. As result, we
obtain that the series $\sum\nolimits_{i\in \mathbf{Z}}\mu _{i+1}^{\ast
}e^{-S_{i}^{\ast }}$ converges $\mathbf{P}^{\ast }$-a.s. Thus, $\Sigma
_{1}<+\infty $ $\mathbf{P}^{\ast }$-a.s.

Fix $i\in \mathbf{Z}$. If the random environment $\mathcal{E}^{\ast }$ is
fixed, the random sequence $\left\{ \eta _{i}^{\ast };\,a_{i,j}^{\ast
}Z_{i,j}^{\ast },\,j\in \mathbf{N}_{i+1}\right\} $ is a martingale.
Therefore
\begin{equation}
\mathbf{E}^{\ast }\left( \left. a_{i,j}^{\ast }Z_{i,j}^{\ast }\,%
\right\vert \,\mathcal{E}^{\ast }\right) =\mu _{i+1}^{\ast }
\label{48}
\end{equation}%
for $j\in \mathbf{N}_{i+1}$. By (5), (48) using Fatou's lemma we obtain that%
\[
\mathbf{E}^{\ast }\left( \left. \zeta _{i}^{\ast }\,\right\vert\, \mathcal{E}^{\ast }\right) \leq \liminf_{j\rightarrow \infty }\mathbf{E}%
^{\ast }\left( \left. a_{i,j}^{\ast }Z_{i,j}^{\ast }\,\right\vert
\,\mathcal{E}^{\ast }\right) =\mu _{i+1}^{\ast }
\]%
and, consequently,
\begin{equation}
\mathbf{E}^{\ast }\left( \left. \zeta _{i}^{\ast }e^{-S_{i}^{\ast }}\,%
\right\vert \,\mathcal{E}^{\ast }\right) =e^{-S_{i}^{\ast }}\mathbf{E}%
^{\ast }\left( \left. \zeta _{i}^{\ast }\,\right\vert \,\mathcal{%
E}^{\ast }\right) \leq \mu _{i+1}^{\ast }e^{-S_{i}^{\ast }}.  \label{49}
\end{equation}%
We have proved that the series $\sum\nolimits_{i\in \mathbf{Z}}\mu
_{i+1}^{\ast }e^{-S_{i}^{\ast }}$ converges $\mathbf{P}^{\ast }$-a.s. This
fact combined with (49) implies convergence of the series $%
\sum\nolimits_{i\in \mathbf{Z}}\mathbf{E}^{\ast }\left( \left. \zeta
_{i+1}^{\ast }e^{-S_{i}^{\ast }}\,\right\vert \,\mathcal{E}%
^{\ast }\right) $ $\mathbf{P}^{\ast }$-a.s. Since the random variables $%
\zeta _{i+1}^{\ast }e^{-S_{i}^{\ast }}$ are nonnegative, it follows that the
series $\sum\nolimits_{i\in \mathbf{Z}}\zeta _{i}^{\ast }e^{-S_{i}^{\ast }}$
converges a.s. for any fixed environment $\mathcal{E}^{\ast }$. Hence, $%
\Sigma _{2}<+\infty $ $\mathbf{P}^{\ast }$-a.s.

The lemma is proved.\medskip

Set%
\[
\Sigma _{1}^{\left( 1\right) }=\sum\limits_{i=0}^{\infty }\mu
_{i+1}^{+}e^{-S_{i}^{+}}=\sum\limits_{i\in \mathbf{N}_{0}}\mu _{i+1}^{\ast
}e^{-S_{i}^{\ast }},
\]%
\[
\Sigma _{1}^{\left( 2\right) }=\sum\limits_{i=1}^{\infty }\mu
_{i}^{-}e^{S_{i}^{-}}=\sum\limits_{i\in \mathbf{Z\setminus N}_{0}}\mu
_{i+1}^{\ast }e^{-S_{i}^{\ast }}.
\]%
Clearly,%
\begin{equation}
\Sigma _{1}=\Sigma _{1}^{\left( 1\right) }+\Sigma _{1}^{\left( 2\right) }
\label{50}
\end{equation}%
and by virtue of Lemma 4 $\mathbf{P}^{\ast }$-a.s.
\begin{equation}
\Sigma _{1}^{\left( 1\right) }<+\infty ,\qquad\Sigma _{1}^{\left(
2\right) }<+\infty .  \label{51}
\end{equation}%
\medskip

\textbf{Lemma 5}. \textit{If the conditions} \textit{of Theorem 1 are}
\textit{satisfied}, \textit{then} $\mathbf{P}^{\ast }$-\textit{a.s.},
\textit{as }$n\rightarrow \infty $,%
\begin{equation}
\left\{ \left. \sum\limits_{i=0}^{n-1}\mu _{i+1}e^{-S_{i}}\,%
\right\vert \,L_{n}\,\geq 0\right\} \stackrel{D}{\to}%
\Sigma _{1}^{\left( 1\right)
},  \label{52}
\end{equation}%
\begin{equation}
\left\{ \left. \sum\limits_{i=1}^{n-1}\mu _{i}e^{S_{i}}\,\right\vert
\,M_{n}\,<0\right\} \stackrel{D}{\to}\Sigma _{1}^{\left( 2\right) }.
\label{53}
\end{equation}%
\medskip

\textit{Proof}. Let $f$ $:$\ $\mathbf{R}\rightarrow \mathbf{R}$ be a bounded
and continuous function. By virtue of (3) for fixed $k\in \mathbf{N}$%
\[
\left\{ \left. f\left( \sum\limits_{i=0}^{k}\mu _{i+1}e^{-S_{i}}\,%
\right) \,\right\vert \,L_{n}\,\geq 0\right\} \stackrel{D}{\to}%
f\left( \sum\limits_{i=0}^{k}\mu _{i+1}^{+}e^{-S_{i}^{+}}\right)
\]%
as $n\rightarrow \infty $. Recalling (51) we conclude that
\[
{\lim_{k\rightarrow \infty }}f\left( \sum\limits_{i=0}^{k}\mu
_{i+1}^{+}e^{-S_{i}^{+}}\right) =f\left( \Sigma _{1}^{\left(
1\right) }\right)
\]%
$\mathbf{P}^{\ast }$-a.s. From these two facts, in view of Lemma 2.5 of \cite%
{Afanasyev2005645}, it follows that%
\[
\left\{ \left. f\left( \sum\limits_{i=0}^{n-1}\mu _{i+1}e^{-S_{i}}\,%
\right) \,\right\vert \,L_{n}\,\geq 0\right\} \stackrel{D}{\to}%
f\left( \Sigma _{1}^{\left( 1\right) }\right) .
\]%
Thus, relation (52) is true. Relation (53) can be proved by similar
arguments.

The lemma is proved.\medskip

\textbf{Remark 4}. It is not difficult to verify that if we combine the
left-hand sides of relations (3) and (52) (or (4) and (53)), then the
respective statements concerning convergence in distribution of the four
dimensional tuples of the random elements given $L_{n}$ $\geq 0$ (or $M_{n}$
$<0$) are still force.\medskip

Set for $n\in \mathbf{N}$%
\[
\mu _{i,n}^{\prime }=\left\{
\begin{array}{c}
\mu _{\tau _{n}+i},\qquad i\in \mathbf{N}_{\left( -\tau _{n}\right) },
\\
0,\qquad i\in \mathbf{Z\setminus N}_{\left( -\tau _{n}\right) }.%
\end{array}%
\right.
\]%
Let%
\[
\Sigma _{1}^{\left( 1\right) }\left( n\right) =\sum\limits_{j=0}^{n-1-\tau
_{n}}\mu _{j+1,n}^{\prime }e^{-S_{j,n}^{\prime }},\qquad\Sigma
_{1}^{\left( 2\right) }\left( n\right) =\sum\limits_{j=1}^{\tau _{n}}\mu
_{-j+1,n}^{\prime }e^{-S_{j,n}^{\prime }}.
\]%
\medskip

\textbf{Lemma 6}. \textit{If the conditions} \textit{of Theorem 1 are}
\textit{satisfied}, \textit{then} $\mathbf{P}^{\ast }$-\textit{a.s.},
\textit{as }$n\rightarrow \infty $,%
\begin{equation}
\left( \left\{ \left( \mu _{i,n}^{\prime },S_{i,n}^{\prime }\right) ,\,%
i\in \mathbf{N}_{0}\right\} ,\,\Sigma _{1}^{\left( 1\right) }\left(
n\right) \right) \stackrel{D}{\to}\left( \left\{ \left( \mu
_{i}^{\ast },S_{i}^{\ast }\right) ,\,i\in \mathbf{N}_{0}\right\} ,%
\,\Sigma _{1}^{\left( 1\right) }\right) ,  \label{54}
\end{equation}%
\begin{equation}
\left( \left\{ \left( \mu _{-i,n}^{\prime },S_{-i,n}^{\prime
}\right) , i\in \mathbf{N}\right\} ,\,\Sigma _{1}^{\left( 2\right)
}\left( n\right) \right) \stackrel{D}{\to}\left( \left\{ \left( \mu
_{-i}^{\ast },S_{-i}^{\ast }\right) ,\,i\in \mathbf{N}\right\} ,\,
\Sigma _{1}^{\left( 2\right) }\right) . \label{55}
\end{equation}%
\textit{Moreover, the left-hand sides of these relations are asymptotically
independent\medskip .\ }

\textit{Proof}. We prove for simplicity only convergence in distribution
(for a fixed $i$) of the random sequences $\left( \mu _{i,n}^{\prime
},S_{i,n}^{\prime },\Sigma _{1}^{\left( 1\right) }\left( n\right) \right) \ $%
and $\left( \mu _{-i,n}^{\prime },S_{-i,n}^{\prime },\Sigma _{1}^{\left(
2\right) }\left( n\right) \right) $, as $n\rightarrow \infty $.

Fix $i\in \mathbf{N}_{0}$. Similarly to relation (7), we can show that, for
any bounded and continuous function $f$ $:\mathbf{R}^{3}\rightarrow \mathbf{R%
}$,%
\begin{eqnarray}
&&\mathbf{E}\left[ f\left( \mu _{i,n}^{\prime },S_{i,n}^{\prime
},\sum\limits_{j=0}^{n-1-\tau _{n}}\mu _{j+1,n}^{\prime }e^{-S_{j,n}^{\prime
}}\right) ;\,\tau _{n}+i\leq n\right] \nonumber\\
&=&\sum\limits_{k=0}^{n-i}\mathbf{E}\left( \left. f\left( \mu
_{i},S_{i},\sum\limits_{j=0}^{n-1-k}\mu _{j+1}e^{-S_{j}}\right) \,%
\right\vert \,L_{n-k}\,\geq 0\right) \mathbf{P}\left( \tau
_{n}=k\right)  \label{56}
\end{eqnarray}%
for $n\geq i$. Repeating the arguments of Lemma 1 and using Lemma 5 and
Remark 4, we can deduce from (56) that%
\begin{eqnarray*}
&&{\lim_{n\rightarrow \infty } }\mathbf{E}f\left( \mu _{i,n}^{\prime
},S_{i,n}^{\prime },\sum\limits_{j=0}^{n-1-\tau _{n}}\mu
_{j+1,n}^{\prime
}e^{-S_{j,n}^{\prime }}\right) \\
&=&\mathbf{E}^{\ast }f\left( \mu _{i}^{+},S_{i}^{+},\sum\limits_{j\in
\mathbf{N}_{0}}\mu _{j+1}^{+}e^{-S_{j}^{+}}\right) =\mathbf{E}^{\ast
}f\left( \mu _{i}^{\ast },S_{i}^{\ast },\Sigma _{1}^{\left( 1\right)
}\right) .
\end{eqnarray*}%
Thus, relation (54) is proved.

Now fix $i\in \mathbf{N}$. It is easy to show (see the proof of relation
(13)) that for $n\geq i$%
\begin{eqnarray*}
&&\mathbf{E}\left[ f\left( \mu _{-i,n}^{\prime },S_{-i,n}^{\prime
},\sum\limits_{j=1}^{\tau _{n}}\mu _{-j+1,n}^{\prime }e^{-S_{-j,n}^{\prime
}}\right) ;\,\tau _{n}-i\geq 0\right] \\
&=&\sum\limits_{k=i}^{n}\mathbf{E}\left( \left. f\left( \mu
_{i+1},-S_{i},\sum\limits_{j=1}^{k}\mu _{j}e^{S_{j}}\right) \,%
\right\vert \,M_{k}\,<0\right) \mathbf{P}\left( \tau
_{n}=k\right)
\end{eqnarray*}%
and therefore (see Lemma 5 and Remark 4)%
\begin{eqnarray*}
&&{\lim_{n\rightarrow \infty } }\mathbf{E}f\left( \mu
_{-i,n}^{\prime },S_{-i,n}^{\prime },\sum\limits_{j=1}^{\tau
_{n}}\mu
_{-j+1,n}^{\prime }e^{-S_{-j,n}^{\prime }}\right) \\
&=&\mathbf{E}^{\ast }f\left( \mu
_{i+1}^{-},-S_{i}^{-},\sum\limits_{j=1}^{\infty }\mu
_{j}^{-}e^{S_{j}^{-}}\right) =\mathbf{E}^{\ast }f\left( \mu _{-i}^{\ast
},S_{-i}^{\ast },\Sigma _{1}^{\left( 2\right) }\right) .
\end{eqnarray*}%
This proves (55). The asymptotic independence of the left-hand sides of
relations (54) and (55) is obvious.

The lemma is proved.\medskip

\textbf{Remark 5}. It is not difficult to verify that statement (54) admits
the following generalization: for any $a\leq 0$ and $b>0$, as $n\rightarrow
\infty $,%
\begin{eqnarray*}
&&\left( \left\{ \left( \mu _{i,n}^{\prime },S_{i,n}^{\prime }\right) , i\in \mathbf{N}_{0}\right\} ,\,\Sigma _{1}^{\left( 1\right) }\left(
n\right) \,\left\vert \,\frac{L_{n}}{C_{n}}\leq a,\,\frac{%
S_{n}-L_{n}}{C_{n}}\leq b\right. \right) \\
&&\stackrel{D}{\to}\left( \left\{ \left( \mu _{i}^{\ast
},S_{i}^{\ast }\right) ,\,i\in \mathbf{N}_{0}\right\} ,\,\Sigma
_{1}^{\left( 1\right) }\right) .
\end{eqnarray*}%
Statement (55) allows for a similar generalization.\medskip

\textbf{Lemma 7}. \textit{If the conditions} \textit{of Theorem 1 are}
\textit{satisfied}, \textit{then}, \textit{as }$n\rightarrow \infty $,%
\[
\left\{ \frac{b_{n}-b_{\tau _{n}+i}}{b_{n}},\,i\in \mathbf{N}%
_{0}\right\} \stackrel{D}{\to}\left\{ \frac{\sum\nolimits_{j=i}^{%
\infty }\mu _{j+1}^{+}\exp \left( -S_{j}^{+}\right) }{\Sigma _{1}},\,%
i\in \mathbf{N}_{0}\right\} ,
\]%
\[
\left\{ \frac{b_{\tau _{n}-i}}{b_{n}},\,i\in \mathbf{N}\right\}
\stackrel{D}{\to}\left\{ \frac{\sum\nolimits_{j=i+1}^{\infty }\mu
_{j}^{-}\exp \left( S_{j}^{-}\right) }{\Sigma _{1}},\,i\in \mathbf{N}%
\right\} ,
\]%
\[
\left\{ \frac{a_{\tau _{n}+i}}{b_{n}},\,i\in \mathbf{N}_{0}\right\}
\stackrel{D}{\to}\left\{ \frac{\exp \left( -S_{i}^{+}\right) }{%
\Sigma _{1}},\,i\in \mathbf{N}_{0}\right\} ,
\]%
\[
\left\{ \frac{a_{\tau _{n}-i}}{b_{n}},\,i\in \mathbf{N}\right\}
\stackrel{D}{\to}\left\{ \frac{\exp \left( S_{i}^{-}\right) }{\Sigma
_{1}},\,i\in \mathbf{N}\right\} .
\]%
\medskip

\textit{Proof}. To simplify the presentation we check the first statement
only. Moreover, we prove only convergence of one-dimensional distributions.
Fix $i\in \mathbf{N}_{0}$. Note that for $\tau _{n}+i\leq n$%
\begin{eqnarray*}
\frac{b_{\tau _{n}+i}}{b_{n}} &=&\frac{\sum\nolimits_{j=0}^{\tau
_{n}+i-1}\mu _{j+1}\exp \left( -S_{j}\right) }{\sum\nolimits_{j=0}^{n-1}\mu
_{j+1}\exp \left( -S_{j}\right) }=\frac{\sum\nolimits_{j=0}^{\tau
_{n}+i-1}\mu _{j+1}\exp \left( -\left( S_{j}-S_{\tau _{n}}\right) \right) }{%
\sum\nolimits_{j=0}^{n-1}\mu _{j+1}\exp \left( -\left( S_{j}-S_{\tau
_{n}}\right) \right) } \\
&=&\frac{\sum\nolimits_{j=0}^{i-1}\mu _{j+1,n}^{\prime }\exp \left(
-S_{j,n}^{\prime }\right) +\Sigma _{1}^{\left( 2\right) }\left( n\right) }{%
\Sigma _{1}^{\left( 1\right) }\left( n\right) +\Sigma _{1}^{\left( 2\right)
}\left( n\right) }.
\end{eqnarray*}%
Since the last expression is a bounded continuous function of the random
element mentioned in Lemma 6, it follows that
\[
\frac{b_{\tau _{n}+i}}{b_{n}}\stackrel{D}{\to}\frac{%
\sum\nolimits_{j=0}^{i-1}\mu _{j+1}^{+}\exp \left( -S_{j}^{+}\right) +\Sigma
_{1}^{\left( 2\right) }}{\Sigma _{1}^{\left( 1\right) }+\Sigma _{1}^{\left(
2\right) }}
\]%
as $n\rightarrow \infty $. Whence, taking into account (50) we obtain the
required relation. The remaining three statements may be proved by similar
arguments.

The lemma is proved.\medskip

\textbf{Remark 6}. We can construct a new random element by combining the
left-hand sides of all the relations included in Lemmas 3 and 7. It is not
difficult to prove convergence in distribution of the sequence of these
random elements to a random element constructed by the right-hand sides of
the corresponding relations of Lemmas 3 and 7. Moreover, a random element
constructed by the left-hand sides is asymptotically independent, as $%
n\rightarrow \infty $, of the random event
\[
\left\{ C_{n}^{-1}L_{n}\leq a,\,C_{n}^{-1}\left( S_{n}-L_{n}\right)
\leq b\right\}
\]%
for any $a\leq 0\ $and $b>0$.\medskip

\textbf{3}. \textbf{Proof of the main result\medskip }

\textit{First part}. We establish convergence of one-dimensional
distributions: if $t>0$, then, as $n\rightarrow \infty $,%
\begin{equation}
\frac{a_{\left\lfloor nt\right\rfloor }}{b_{\left\lfloor nt\right\rfloor }}%
Z_{\left\lfloor nt\right\rfloor }\stackrel{D}{\to}\frac{\Sigma _{2}}{%
\Sigma _{1}}.  \label{57}
\end{equation}%
Set for $r\in \mathbf{N}$%
\[
U_{r}^{\left( i\right) }=\sum\limits_{j=\tau _{r}-i}^{\tau _{r}+i-1}Z_{j,r},
\]%
\[
V_{r}^{\left( i\right) }=\sum\limits_{j=0}^{\tau
_{r}-i-1}Z_{j,r}+\sum\limits_{j=\tau _{r}+i}^{r-1}Z_{j,r}.
\]%
It is clear that for $i\in \mathbf{N}$%
\begin{equation}
Z_{\left\lfloor nt\right\rfloor }=\sum\limits_{j=0}^{\left\lfloor
nt\right\rfloor -1}Z_{j,\left\lfloor nt\right\rfloor }=U_{\left\lfloor
nt\right\rfloor }^{\left( i\right) }+V_{\left\lfloor nt\right\rfloor
}^{\left( i\right) }.  \label{58}
\end{equation}

Note that
\begin{equation}
\mathbf{E}\left( \left. a_{j,\left\lfloor nt\right\rfloor }Z_{j,\left\lfloor
nt\right\rfloor }\,\right\vert \,Q_{1,\left\lfloor
nt\right\rfloor }\right) =\mu _{j+1},  \label{59}
\end{equation}%
if $1\leq $ $j<\left\lfloor nt\right\rfloor $. Observing that $%
a_{\left\lfloor nt\right\rfloor }=a_{j}a_{j,\left\lfloor nt\right\rfloor }$
for $1\leq $ $j<\left\lfloor nt\right\rfloor $ we obtain by (59) that%
\begin{eqnarray}
&&\mathbf{E}\left( \frac{a_{\left\lfloor nt\right\rfloor }}{b_{\left\lfloor
nt\right\rfloor }}V_{\left\lfloor nt\right\rfloor }^{\left( i\right) }\right)
\nonumber \\
&=&\mathbf{E}b_{\left\lfloor nt\right\rfloor }^{-1}\left(
\sum\limits_{j=0}^{\tau _{\left\lfloor nt\right\rfloor
}-i-1}a_{j}a_{j,\left\lfloor nt\right\rfloor }Z_{j,\left\lfloor
nt\right\rfloor }+\sum\limits_{j=\tau _{\left\lfloor nt\right\rfloor
}+i}^{\left\lfloor nt\right\rfloor -1}a_{j}a_{j,\left\lfloor nt\right\rfloor
}Z_{j,\left\lfloor nt\right\rfloor }\right) \nonumber\\
&=&\mathbf{E}b_{\left\lfloor nt\right\rfloor }^{-1}\left(
\sum\limits_{j=0}^{\tau _{\left\lfloor nt\right\rfloor }-i-1}\mu
_{j+1}a_{j}+\sum\limits_{j=\tau _{\left\lfloor nt\right\rfloor
}+i}^{\left\lfloor nt\right\rfloor -1}\mu _{j+1}a_{j}\right) \nonumber\\
&=&\mathbf{E}\frac{b_{\tau _{\left\lfloor nt\right\rfloor }-i}+\left(
b_{\left\lfloor nt\right\rfloor }-b_{\tau _{\left\lfloor nt\right\rfloor
}+i}\right) }{b_{\left\lfloor nt\right\rfloor }}.  \label{60}
\end{eqnarray}%
Applying Lemma 7 to the right-hand side of (60), we conclude that%
\[
\lim_{n\rightarrow \infty }\mathbf{E}\left( \frac{a_{\left\lfloor
nt\right\rfloor }}{b_{\left\lfloor nt\right\rfloor }}V_{\left\lfloor
nt\right\rfloor }^{\left( i\right) }\right) =\frac{\sum\nolimits_{j=i}^{%
\infty }\mu _{j+1}^{+}\exp \left( -S_{j}^{+}\right)
+\sum\nolimits_{j=i+1}^{\infty }\mu _{j}^{-}\exp \left( S_{j}^{-}\right) }{%
\Sigma _{1}}
\]%
and, therefore (see Lemma 4),%
\begin{equation}
\lim_{i\rightarrow \infty }\lim_{n\rightarrow \infty }\mathbf{E}\left( \frac{%
a_{\left\lfloor nt\right\rfloor }}{b_{\left\lfloor nt\right\rfloor }}%
V_{\left\lfloor nt\right\rfloor }^{\left( i\right) }\right) =0.  \label{61}
\end{equation}%
By Markov inequality for any $\varepsilon >0$%
\[
\mathbf{P}\left( \frac{a_{\left\lfloor nt\right\rfloor }}{b_{\left\lfloor
nt\right\rfloor }}V_{\left\lfloor nt\right\rfloor }^{\left( i\right) }\geq
\varepsilon \right) \leq \varepsilon ^{-1}\mathbf{E}\left( \frac{%
a_{\left\lfloor nt\right\rfloor }}{b_{\left\lfloor nt\right\rfloor }}%
V_{\left\lfloor nt\right\rfloor }^{\left( i\right) }\right) .
\]%
Hence, taking into account (61) we obtain that%
\begin{equation}
\lim_{i\rightarrow \infty }\lim_{n\rightarrow \infty }\mathbf{P}\left( \frac{%
a_{\left\lfloor nt\right\rfloor }}{b_{\left\lfloor nt\right\rfloor }}%
V_{\left\lfloor nt\right\rfloor }^{\left( i\right) }\geq \varepsilon \right)
=0.  \label{62}
\end{equation}

Observe that we may assume in the sequel that $i\leq \tau _{\left\lfloor
nt\right\rfloor }<\left\lfloor nt\right\rfloor -i$ (see the proof of Lemma
1). Note that
\[
U_{\left\lfloor nt\right\rfloor }^{\left( i\right)
}=\sum\limits_{j=-i}^{i-1}Z_{\tau _{\left\lfloor nt\right\rfloor
}+j,\left\lfloor nt\right\rfloor }=\sum\limits_{j=-i}^{i-1}Z_{j,\left\lfloor
nt\right\rfloor }^{\prime }
\]%
and, therefore,%
\begin{equation}
\frac{a_{\left\lfloor nt\right\rfloor }}{b_{\left\lfloor nt\right\rfloor }}%
U_{\left\lfloor nt\right\rfloor }^{\left( i\right) }=\sum\limits_{j=-i}^{i-1}%
\frac{a_{\tau _{\left\lfloor nt\right\rfloor }+j}}{b_{\left\lfloor
nt\right\rfloor }}a_{j,\left\lfloor nt\right\rfloor }^{\prime
}Z_{j,\left\lfloor nt\right\rfloor }^{\prime }.  \label{63}
\end{equation}%
Applying Lemmas 3, 7 and Remark 6 to relation (63), we obtain that, as $%
n\rightarrow \infty $,%
\begin{equation}
\frac{a_{\left\lfloor nt\right\rfloor }}{b_{\left\lfloor nt\right\rfloor }}%
U_{\left\lfloor nt\right\rfloor }^{\left( i\right) }\stackrel{D}{\to}%
\frac{1}{\Sigma _{1}}\sum\limits_{j=-i}^{i-1}\zeta _{j}^{\ast
}e^{-S_{j}^{\ast }}.  \label{64}
\end{equation}%
Hence, for all but a countable set of $x\geq 0$%
\begin{equation}
\lim_{n\rightarrow \infty }\mathbf{P}\left( \frac{a_{\left\lfloor
nt\right\rfloor }}{b_{\left\lfloor nt\right\rfloor }}U_{\left\lfloor
nt\right\rfloor }^{\left( i\right) }\leq x\right) =\mathbf{P}\left( \frac{1}{%
\Sigma _{1}}\sum\limits_{j=-i}^{i-1}\zeta _{j}^{\ast }e^{-S_{j}^{\ast }}\leq
x\right) .  \label{65}
\end{equation}%
In view of Lemma 4%
\begin{equation}
\lim_{i\rightarrow \infty }\mathbf{P}\left( \frac{1}{\Sigma _{1}}%
\sum\limits_{j=-i}^{i-1}\zeta _{i}^{\ast }e^{-S_{i}^{\ast }}\leq x\right) =%
\mathbf{P}\left( \frac{\Sigma _{2}}{\Sigma _{1}}\leq x\right) .  \label{66}
\end{equation}%
We obtain by (65) and (66) that%
\begin{equation}
\lim_{i\rightarrow \infty }\lim_{n\rightarrow \infty }\mathbf{P}\left( \frac{%
a_{\left\lfloor nt\right\rfloor }}{b_{\left\lfloor nt\right\rfloor }}%
U_{\left\lfloor nt\right\rfloor }^{\left( i\right) }\leq x\right) =\mathbf{P}%
\left( \frac{\Sigma _{2}}{\Sigma _{1}}\leq x\right) .  \label{67}
\end{equation}%
It follows from (58), (62) and (67) that for all but a countable set of $%
x\geq 0$%
\[
\lim_{n\rightarrow \infty }\mathbf{P}\left( \frac{a_{\left\lfloor
nt\right\rfloor }}{b_{\left\lfloor nt\right\rfloor }}Z_{\left\lfloor
nt\right\rfloor }\leq x\right) =\mathbf{P}\left( \frac{\Sigma _{2}}{\Sigma
_{1}}\leq x\right) .
\]%
This proves (57).\medskip

\textbf{Remark 7}. It is not difficult to verify that relation (64) admits
the following generalization: for any $a\leq 0$ and $b>0$, as $n\rightarrow
\infty $,%
\[
\left\{ \left. \frac{a_{\left\lfloor nt\right\rfloor }}{b_{\left\lfloor
nt\right\rfloor }}U_{\left\lfloor nt\right\rfloor }^{\left( i\right) }\,\right\vert \,\frac{L_{\left\lfloor nt\right\rfloor }}{C_{n}}\leq a,%
\,\frac{S_{\left\lfloor nt\right\rfloor }-L_{\left\lfloor
nt\right\rfloor }}{C_{n}}\leq b\right\} \stackrel{D}{\to}\frac{1}{%
\Sigma _{1}}\sum\limits_{j=-i}^{i+1}\zeta _{j}^{\ast }e^{-S_{j}^{\ast }}.
\]%
\medskip

\textit{Second part}. Now we establish convergence of two-dimensional
distributions. Select $0<t_{1}<t_{2}$, fix an $\varepsilon >0$ and introduce
the following random events:%
\[
A_{n,\varepsilon }=\left\{ L_{\left\lfloor nt_{1}\right\rfloor
}>L_{\left\lfloor nt_{1}\right\rfloor ,\left\lfloor nt_{2}\right\rfloor
}+\varepsilon C_{n}\right\} ,
\]%
\[
B_{n,\varepsilon }=\left\{ L_{\left\lfloor nt_{1}\right\rfloor
}<L_{\left\lfloor nt_{1}\right\rfloor ,\left\lfloor nt_{2}\right\rfloor
}-\varepsilon C_{n}\right\} ,
\]%
\[
D_{n,\varepsilon }=\left\{ \left\vert L_{\left\lfloor nt_{1}\right\rfloor
}-L_{\left\lfloor nt_{1}\right\rfloor ,\left\lfloor nt_{2}\right\rfloor
}\right\vert \leq \varepsilon C_{n}\right\} .
\]%
We show that, as $n\rightarrow \infty $,%
\begin{equation}
\left\{ \left. \frac{a_{\left\lfloor nt_{1}\right\rfloor }}{b_{\left\lfloor
nt_{1}\right\rfloor }}Z_{\left\lfloor nt_{1}\right\rfloor },\frac{%
a_{\left\lfloor nt_{2}\right\rfloor }}{b_{\left\lfloor nt_{2}\right\rfloor }}%
Z_{\left\lfloor nt_{2}\right\rfloor }\,\right\vert \,%
A_{n,\varepsilon }\right\} \stackrel{D}{\to}\left( \gamma
_{1},\gamma _{2}\right) ,  \label{68}
\end{equation}%
\begin{equation}
\left\{ \left. \frac{a_{\left\lfloor nt_{1}\right\rfloor }}{b_{\left\lfloor
nt_{1}\right\rfloor }}Z_{\left\lfloor nt_{1}\right\rfloor },\frac{%
a_{\left\lfloor nt_{2}\right\rfloor }}{b_{\left\lfloor nt_{2}\right\rfloor }}%
Z_{\left\lfloor nt_{2}\right\rfloor }\,\right\vert \,%
B_{n,\varepsilon }\right\} \stackrel{D}{\to}\left( \gamma
_{1},\gamma _{1}\right) ,  \label{69}
\end{equation}%
where $\gamma _{1},\gamma _{2}$ are independent random variables and
$\gamma _{1}\stackrel{D}{=}\gamma _{2}\stackrel{D}{=}\Sigma
_{2}/\Sigma _{1}$.

First we establish (68). To this aim we prove that, for any fixed $i\in
\mathbf{N}$ and for all but a countable set of $\left( x_{1},x_{2}\right) $
with $x_{1},x_{2}\geq 0$,%
\begin{eqnarray}
&&\lim_{n\rightarrow \infty }\mathbf{P}\left( \left. \frac{a_{\left\lfloor
nt_{1}\right\rfloor }}{b_{\left\lfloor nt_{1}\right\rfloor }}U_{\left\lfloor
nt_{1}\right\rfloor }^{\left( i\right) }\leq x_{1},\,\frac{%
b_{\left\lfloor nt_{2}\right\rfloor }}{a_{\left\lfloor nt_{2}\right\rfloor }}%
U_{\left\lfloor nt_{2}\right\rfloor }^{\left( i\right) }\leq x_{2}\,%
\right\vert \,A_{n,\varepsilon }\right) \nonumber\\
&=&\mathbf{P}\left( \frac{1}{\Sigma _{1}}\sum\limits_{j=-i}^{i+1}\zeta
_{j}^{\ast }e^{-S_{j}^{\ast }}\leq x_{1}\right) \mathbf{P}\left( \frac{1}{%
\Sigma _{1}}\sum\limits_{j=-i}^{i+1}\zeta _{j}^{\ast }e^{-S_{j}^{\ast }}\leq
x_{2}\right) .  \label{70}
\end{eqnarray}

Provided the random event $A_{n,\varepsilon }$ occurred, it follows that, as
$n\rightarrow \infty $,%
\[
\frac{b_{\left\lfloor nt_{2}\right\rfloor }}{a_{\left\lfloor
nt_{2}\right\rfloor }}\sim \frac{b_{\left\lfloor nt_{2}\right\rfloor
}-b_{\left\lfloor nt_{1}\right\rfloor }}{a_{\left\lfloor nt_{2}\right\rfloor
}}=\frac{\widetilde{b}_{\left\lfloor nt_{2}\right\rfloor -\left\lfloor
nt_{1}\right\rfloor }}{\widetilde{a}_{\left\lfloor nt_{2}\right\rfloor
-\left\lfloor nt_{1}\right\rfloor }},
\]%
where the values $\widetilde{a}_{\left\lfloor nt_{2}\right\rfloor
-\left\lfloor nt_{1}\right\rfloor }\ $and $\widetilde{b}_{\left\lfloor
nt_{2}\right\rfloor -\left\lfloor nt_{1}\right\rfloor }$ are constructed by
the random environment $\widetilde{Q}_{i}:=Q_{\left\lfloor
nt_{1}\right\rfloor +i}$, $i=1,\ldots ,\left\lfloor nt_{2}\right\rfloor
-\left\lfloor nt_{1}\right\rfloor $, just as the values $a_{\left\lfloor
nt_{2}\right\rfloor -\left\lfloor nt_{1}\right\rfloor }$ and $%
b_{\left\lfloor nt_{2}\right\rfloor -\left\lfloor nt_{1}\right\rfloor }$ are
constructed by the random environment $Q_{1,\left\lfloor nt_{2}\right\rfloor
-\left\lfloor nt_{1}\right\rfloor }$.

Further, given $A_{n,\varepsilon }$, the inequality $\tau _{\left\lfloor
nt_{2}\right\rfloor }>\tau _{\left\lfloor nt_{1}\right\rfloor }$ is true (we
may assume that $\left\lfloor nt_{2}\right\rfloor -i>\tau _{\left\lfloor
nt_{2}\right\rfloor }>\tau _{\left\lfloor nt_{1}\right\rfloor }+i$). Thus,
if the random environment $\left\{ Q_{n},\,n\in \mathbf{N}\right\} $
is fixed, the distribution of the random variable $U_{\left\lfloor
nt_{1}\right\rfloor }^{\left( i\right) }$ is completely determined by the
random environment $Q_{1,\left\lfloor nt_{1}\right\rfloor }$ and the
distribution of the random variable $U_{\left\lfloor nt_{2}\right\rfloor
}^{\left( i\right) }$ is completely determined by the random environment $%
Q_{\left\lfloor nt_{1}\right\rfloor +1,\left\lfloor nt_{2}\right\rfloor }$.
Moreover, $U_{\left\lfloor nt_{2}\right\rfloor }^{\left( i\right) }=%
\widetilde{U}_{\left\lfloor nt_{2}\right\rfloor -\left\lfloor
nt_{1}\right\rfloor }^{\left( i\right) }$, where $\widetilde{U}%
_{\left\lfloor nt_{2}\right\rfloor -\left\lfloor nt_{1}\right\rfloor
}^{\left( i\right) }$ has the same meaning for the environment $\widetilde{Q}%
_{i}$, $i=1,\ldots ,\left\lfloor nt_{2}\right\rfloor -\left\lfloor
nt_{1}\right\rfloor $, as $U_{\left\lfloor nt_{2}\right\rfloor -\left\lfloor
nt_{1}\right\rfloor }^{\left( i\right) }$ has for the environment $%
Q_{1,\left\lfloor nt_{2}\right\rfloor -\left\lfloor nt_{1}\right\rfloor }$.

Summarizing the arguments above, we see that to prove (70) it is sufficient
to show that%
\begin{eqnarray}
&&\lim_{n\rightarrow \infty }\mathbf{P}\left( \left. \frac{a_{\left\lfloor
nt_{1}\right\rfloor }}{b_{\left\lfloor nt_{1}\right\rfloor }}U_{\left\lfloor
nt_{1}\right\rfloor }^{\left( i\right) }\leq x_{1},\,\frac{\widetilde{a%
}_{\left\lfloor nt_{2}\right\rfloor -\left\lfloor nt_{1}\right\rfloor }}{%
\widetilde{b}_{\left\lfloor nt_{2}\right\rfloor -\left\lfloor
nt_{1}\right\rfloor }}\widetilde{U}_{\left\lfloor nt_{2}\right\rfloor
-\left\lfloor nt_{1}\right\rfloor }^{\left( i\right) }\leq x_{2}\,%
\right\vert \,A_{n,\varepsilon }\right) \nonumber\\
&=&\mathbf{P}\left( \frac{1}{\Sigma _{1}}\sum\limits_{j=-i}^{i+1}\zeta
_{j}^{\ast }e^{-S_{j}^{\ast }}\leq x_{1}\right) \mathbf{P}\left( \frac{1}{%
\Sigma _{1}}\sum\limits_{j=-i}^{i+1}\zeta _{j}^{\ast }e^{-S_{j}^{\ast }}\leq
x_{2}\right) .  \label{71}
\end{eqnarray}

Note that%
\begin{eqnarray*}
&&\mathbf{P}\left( \frac{a_{\left\lfloor nt_{1}\right\rfloor }}{%
b_{\left\lfloor nt_{1}\right\rfloor }}U_{\left\lfloor nt_{1}\right\rfloor
}^{\left( i\right) }\leq x_{1},\,\frac{\widetilde{a}_{\left\lfloor
nt_{2}\right\rfloor -\left\lfloor nt_{1}\right\rfloor }}{\widetilde{b}%
_{\left\lfloor nt_{2}\right\rfloor -\left\lfloor nt_{1}\right\rfloor }}%
\widetilde{U}_{\left\lfloor nt_{2}\right\rfloor -\left\lfloor
nt_{1}\right\rfloor }^{\left( i\right) }\leq x_{2},\,A_{n,\varepsilon
}\right) \\
&=&\int\limits_{-\infty }^{0}\int\limits_{0}^{+\infty }\mathbf{P}\left(
\frac{a_{\left\lfloor nt_{1}\right\rfloor }}{b_{\left\lfloor
nt_{1}\right\rfloor }}U_{\left\lfloor nt_{1}\right\rfloor }^{\left( i\right)
}\leq x_{1},\,\frac{L_{\left\lfloor nt_{1}\right\rfloor }}{C_{n}}\in
da,\,\frac{S_{\left\lfloor nt_{1}\right\rfloor }-L_{\left\lfloor
nt_{1}\right\rfloor }}{C_{n}}\in db\right) \\
&&\times \mathbf{P}\left( \frac{a_{\left\lfloor nt_{2}\right\rfloor
-\left\lfloor nt_{1}\right\rfloor }}{b_{\left\lfloor nt_{2}\right\rfloor
-\left\lfloor nt_{1}\right\rfloor }}U_{\left\lfloor nt_{2}\right\rfloor
-\left\lfloor nt_{1}\right\rfloor }^{\left( i\right) }\leq x_{2},\,%
\frac{L_{\left\lfloor nt_{2}\right\rfloor -\left\lfloor nt_{1}\right\rfloor }%
}{C_{n}}<b-a-\varepsilon \right) .
\end{eqnarray*}%
Hence, taking into account Remark 7 we deduce that, as $n\rightarrow \infty $%
,%
\begin{eqnarray*}
&&\mathbf{P}\left( \frac{a_{\left\lfloor nt_{1}\right\rfloor }}{%
b_{\left\lfloor nt_{1}\right\rfloor }}U_{\left\lfloor nt_{1}\right\rfloor
}^{\left( i\right) }\leq x_{1},\,\frac{\widetilde{a}_{\left\lfloor
nt_{2}\right\rfloor -\left\lfloor nt_{1}\right\rfloor }}{\widetilde{b}%
_{\left\lfloor nt_{2}\right\rfloor -\left\lfloor nt_{1}\right\rfloor }}%
U_{\left\lfloor nt_{2}\right\rfloor }^{\left( i\right) }\leq x_{2},\,%
A_{n,\varepsilon }\right) \\
&\sim &\mathbf{P}\left( \frac{1}{\Sigma _{1}}\sum\limits_{j=-i}^{i+1}\zeta
_{j}^{\ast }e^{-S_{j}^{\ast }}\leq x_{1}\right) \mathbf{P}\left( \frac{1}{%
\Sigma _{1}}\sum\limits_{j=-i}^{i+1}\zeta _{j}^{\ast }e^{-S_{j}^{\ast }}\leq
x_{2}\right) \\
&&\times \int\limits_{-\infty }^{0}\int\limits_{0}^{+\infty }\mathbf{P}%
\left( \frac{L_{\left\lfloor nt_{1}\right\rfloor }}{C_{n}}\in da,\,%
\frac{S_{\left\lfloor nt_{1}\right\rfloor }-L_{\left\lfloor
nt_{1}\right\rfloor }}{C_{n}}\in db\right) \\
&&\times \mathbf{P}\left( \frac{L_{\left\lfloor nt_{2}\right\rfloor
-\left\lfloor nt_{1}\right\rfloor }}{C_{n}}<b-a-\varepsilon \right) .
\end{eqnarray*}%
Since the last integral is equal to $\mathbf{P}\left( A_{n,\varepsilon
}\right) $, we obtain (71) and, as result, the required relation (70).

It follows from (70) that (see (67))%
\begin{eqnarray}
&&\lim_{i\rightarrow \infty }\lim_{n\rightarrow \infty }\mathbf{P}\left(
\left. \frac{a_{\left\lfloor nt_{1}\right\rfloor }}{b_{\left\lfloor
nt_{1}\right\rfloor }}U_{\left\lfloor nt_{1}\right\rfloor }^{\left( i\right)
}\leq x_{1},\,\frac{b_{\left\lfloor nt_{2}\right\rfloor }}{%
a_{\left\lfloor nt_{k}\right\rfloor }}U_{\left\lfloor nt_{2}\right\rfloor
}^{\left( i\right) }\leq x_{2}\,\right\vert \,A_{n,\varepsilon
}\right) \nonumber\\
&=&\mathbf{P}\left( \frac{\Sigma _{2}}{\Sigma _{1}}\leq x_{1}\right) \mathbf{%
P}\left( \frac{\Sigma _{2}}{\Sigma _{1}}\leq x_{2}\right) .  \label{72}
\end{eqnarray}

Applying now the same arguments which we have used in First part of the
proof to establish (57) from (67), we obtain (68) from (72).

We now prove (69). To this aim we check that, for any fixed $i\in \mathbf{N}$
and for all but a countable set of $\left( x_{1},x_{2}\right) $ with $%
x_{1},x_{2}\geq 0$,%
\begin{eqnarray}
&&\lim_{n\rightarrow \infty }\mathbf{P}\left( \left. \frac{a_{\left\lfloor
nt_{1}\right\rfloor }}{b_{\left\lfloor nt_{1}\right\rfloor }}U_{\left\lfloor
nt_{1}\right\rfloor }^{\left( i\right) }\leq x_{1},\,\frac{%
b_{\left\lfloor nt_{2}\right\rfloor }}{a_{\left\lfloor nt_{k}\right\rfloor }}%
U_{\left\lfloor nt_{2}\right\rfloor }^{\left( i\right) }\leq x_{2}\,%
\right\vert \,B_{n,\varepsilon }\right) \nonumber\\
&=&\mathbf{P}\left( \frac{1}{\Sigma _{1}}\sum\limits_{j=-i}^{i+1}\zeta
_{j}^{\ast }e^{-S_{j}^{\ast }}\leq \min \left( x_{1},x_{2}\right) \right) .
\label{73}
\end{eqnarray}%
Set%
\[
Z_{i,n}^{\prime }\left( m\right) =Z_{\tau _{n}+i,m},
\]%
\[
U_{n}^{\left( i\right) }\left( m\right)
=\sum\limits_{j=-i}^{i+1}Z_{j,n}^{\prime }\left( m\right) .
\]%
Given that the random event $B_{n,\varepsilon }$ occurred, $\tau
_{\left\lfloor nt_{2}\right\rfloor }=\tau _{\left\lfloor nt_{1}\right\rfloor
}$ and%
\[
\frac{b_{\left\lfloor nt_{2}\right\rfloor }}{a_{\left\lfloor
nt_{2}\right\rfloor }}\sim \frac{b_{\left\lfloor nt_{1}\right\rfloor }}{%
a_{\left\lfloor nt_{2}\right\rfloor }}
\]%
as $n\rightarrow \infty $. Therefore
\[
U_{\left\lfloor nt_{1}\right\rfloor }^{\left( i\right)
}=\sum\limits_{j=-i}^{i+1}Z_{i,\left\lfloor nt_{1}\right\rfloor }^{\prime
}\left( \left\lfloor nt_{1}\right\rfloor \right) =U_{\left\lfloor
nt_{1}\right\rfloor }^{\left( i\right) }\left( \left\lfloor
nt_{1}\right\rfloor \right) ,
\]%
\[
\qquad U_{\left\lfloor nt_{2}\right\rfloor }^{\left( i\right)
}=\sum\limits_{j=-i}^{i+1}Z_{i,\left\lfloor nt_{1}\right\rfloor }^{\prime
}\left( \left\lfloor nt_{2}\right\rfloor \right) =U_{\left\lfloor
nt_{1}\right\rfloor }^{\left( i\right) }\left( \left\lfloor
nt_{2}\right\rfloor \right) .
\]%
Thus, to prove (73) it is sufficient to show that%
\begin{eqnarray}
&&\lim_{n\rightarrow \infty }\mathbf{P}\left( \left. \frac{a_{\left\lfloor
nt_{1}\right\rfloor }}{b_{\left\lfloor nt_{1}\right\rfloor }}U_{\left\lfloor
nt_{1}\right\rfloor }^{\left( i\right) }\left( \left\lfloor
nt_{1}\right\rfloor \right) \leq x_{1},\,\frac{a_{\left\lfloor
nt_{2}\right\rfloor }}{b_{\left\lfloor nt_{1}\right\rfloor }}U_{\left\lfloor
nt_{1}\right\rfloor }^{\left( i\right) }\left( \left\lfloor
nt_{2}\right\rfloor \right) \leq x_{2}\,\right\vert \,%
B_{n,\varepsilon }\right) \nonumber\\
&=&\mathbf{P}\left( \frac{1}{\Sigma _{1}}\sum\limits_{j=-i}^{i+1}\zeta
_{j}^{\ast }e^{-S_{j}^{\ast }}\leq \min \left( x_{1},x_{2}\right) \right) .
\label{74}
\end{eqnarray}%
Applying the arguments similar to those used to establish relation (19), we
can show that%
\[
\left\{ a_{i,m}^{\prime }Z_{i,n}^{\prime }\left( m\right) ,\,i\in
\mathbf{Z}\right\} \stackrel{D}{\to}\left\{ \zeta _{i}^{\ast },\,
i\in \mathbf{Z}\right\} ,
\]%
as $m\geq n\rightarrow \infty $. Moreover,%
\begin{equation}
\left\{ \left( a_{i,n}^{\prime }Z_{i,n}^{\prime }\left( n\right)
,a_{i,m}^{\prime }Z_{i,n}^{\prime }\left( m\right) \right) ,\,i\in
\mathbf{Z}\right\} \stackrel{D}{\to}\left\{ \left( \zeta _{i}^{\ast
},\zeta _{i}^{\ast }\right) ,\,i\in \mathbf{Z}\right\}  \label{75}
\end{equation}%
and the left-hand side of this relation is asymptotically independent from
the random event $\left\{ C_{n}^{-1}L_{n}\leq a,\,C_{n}^{-1}\left(
S_{n}-L_{n}\right) \leq b\right\} $ for any $a\leq 0\ $and $b>0$. It follows
from (75) that (see the proof of (64))%
\begin{equation}
\left( \frac{a_{n}}{b_{n}}U_{n}^{\left( i\right) }\left( n\right) ,\frac{%
a_{m}}{b_{n}}U_{n}^{\left( i\right) }\left( m\right) \right) \stackrel{D}{\to}%
\frac{1}{\Sigma _{1}}\left( \sum\limits_{j=-i}^{i+1}\zeta _{j}^{\ast
}e^{-S_{j}^{\ast }},\sum\limits_{j=-i}^{i+1}\zeta _{j}^{\ast
}e^{-S_{j}^{\ast }}\right) ,  \label{76}
\end{equation}%
as $m\geq n\rightarrow \infty $. From (76) we obtain the desired relation
(74) and, as result, (73). Now statement (69) follows from (73) in a
standard way.

Finally, according to the Skorokhod functional limit theorem (see (1))%
\begin{equation}
\lim_{\varepsilon \rightarrow 0}\lim_{n\rightarrow \infty }\mathbf{P}\left(
A_{n,\varepsilon }\right) =\mathbf{P}\left( L\left( t_{1}\right) >L\left(
t_{1},t_{2}\right) \right) =\mathbf{P}\left( L\left( t_{1}\right) >L\left(
t_{2}\right) \right) ,  \label{77}
\end{equation}%
\begin{equation}
\lim_{\varepsilon \rightarrow 0}\lim_{n\rightarrow \infty }\mathbf{P}\left(
B_{n,\varepsilon }\right) =\mathbf{P}\left( L\left( t_{1}\right) <L\left(
t_{1},t_{2}\right) \right) =\mathbf{P}\left( L\left( t_{1}\right) =L\left(
t_{2}\right) \right) ,  \label{78}
\end{equation}%
where $L\left( t_{1},t_{2}\right) =\inf_{t\in \left[ t_{1},t_{2}\right]
}W\left( t\right) $, and
\begin{equation}
\lim_{\varepsilon \rightarrow 0}\lim_{n\rightarrow \infty }\mathbf{P}\left(
D_{n,\varepsilon }\right) =0.  \label{79}
\end{equation}%
By the total probability formula
\begin{eqnarray}
&&\mathbf{P}\left( \frac{a_{\left\lfloor nt_{1}\right\rfloor }}{%
b_{\left\lfloor nt_{1}\right\rfloor }}Z_{\left\lfloor nt_{1}\right\rfloor
}\leq x_{1},\,\frac{a_{\left\lfloor nt_{2}\right\rfloor }}{%
b_{\left\lfloor nt_{2}\right\rfloor }}Z_{\left\lfloor nt_{2}\right\rfloor
}\leq x_{2}\right) \nonumber\\
&=&\mathbf{P}\left( \left. \frac{a_{\left\lfloor nt_{1}\right\rfloor }}{%
b_{\left\lfloor nt_{1}\right\rfloor }}Z_{\left\lfloor nt_{1}\right\rfloor
}\leq x_{1},\,\frac{a_{\left\lfloor nt_{2}\right\rfloor }}{%
b_{\left\lfloor nt_{2}\right\rfloor }}Z_{\left\lfloor nt_{2}\right\rfloor
}\leq x_{2}\,\right\vert \,A_{n,\varepsilon }\right) \mathbf{P}%
\left( A_{n,\varepsilon }\right) \nonumber\\
&&+\mathbf{P}\left( \left. \frac{a_{\left\lfloor nt_{1}\right\rfloor }}{%
b_{\left\lfloor nt_{1}\right\rfloor }}Z_{\left\lfloor nt_{1}\right\rfloor
}\leq x_{1},\,\frac{a_{\left\lfloor nt_{2}\right\rfloor }}{%
b_{\left\lfloor nt_{2}\right\rfloor }}Z_{\left\lfloor nt_{2}\right\rfloor
}\leq x_{2}\,\right\vert \,B_{n,\varepsilon }\right) \mathbf{P}%
\left( B_{n,\varepsilon }\right) \nonumber\\
&&+\mathbf{P}\left( \left. \frac{a_{\left\lfloor nt_{1}\right\rfloor }}{%
b_{\left\lfloor nt_{1}\right\rfloor }}Z_{\left\lfloor nt_{1}\right\rfloor
}\leq x_{1},\,\frac{a_{\left\lfloor nt_{2}\right\rfloor }}{%
b_{\left\lfloor nt_{2}\right\rfloor }}Z_{\left\lfloor nt_{2}\right\rfloor
}\leq x_{2}\,\right\vert \,D_{n,\varepsilon }\right) \mathbf{P}%
\left( D_{n,\varepsilon }\right) .  \label{80}
\end{eqnarray}%
Combining (68), (69) and (77)-(80) we deduce that%
\begin{eqnarray*}
&&\lim_{n\rightarrow \infty }\mathbf{P}\left( \frac{a_{\left\lfloor
nt_{1}\right\rfloor }}{b_{\left\lfloor nt_{1}\right\rfloor }}Z_{\left\lfloor
nt_{1}\right\rfloor }\leq x_{1},\,\frac{a_{\left\lfloor
nt_{2}\right\rfloor }}{b_{\left\lfloor nt_{2}\right\rfloor }}Z_{\left\lfloor
nt_{2}\right\rfloor }\leq x_{2}\right) \\
&=&\mathbf{P}\left( \gamma _{1}\leq x_{1},\,\gamma _{2}\leq
x_{2}\right) \mathbf{P}\left( L\left( t_{1}\right) >L\left( t_{2}\right)
\right) \\
&&+\mathbf{P}\left( \gamma _{1}\leq x_{1},\,\gamma _{1}\leq
x_{2}\right) \mathbf{P}\left( L\left( t_{1}\right) =L\left( t_{2}\right)
\right) ,
\end{eqnarray*}%
This gives the desired convergence of two-dimensional distributions.\medskip

\textit{Third part}. The proof of convergence of multidimensional
distributions (for dimensions exceeding two) is carried out by induction
using the reasonings of Second part of the proof.\medskip


\bibliography{mybibfile}

\end{document}